\documentclass[11pt,reqno]{amsart}
\usepackage{iftex}
\ifPDFTeX
\usepackage[T1]{fontenc}
\usepackage{lmodern}
\else
\usepackage[no-math]{fontspec}
\fi
\usepackage[a4paper,margin=1in]{geometry}
\usepackage{amsmath,amssymb,amsfonts,amsthm,mathtools,mathrsfs}
\usepackage{enumitem,microtype,xcolor}
\usepackage[numbers,sort&compress]{natbib}
\usepackage[hypertexnames=false]{hyperref}
\allowdisplaybreaks[2]
\hypersetup{bookmarksdepth=2,colorlinks=true,
linkcolor=blue!45!black,citecolor=blue!45!black,urlcolor=blue!45!black,
pdftitle={Deterministic Principal-Ideal Responses of Random Operators},
pdfauthor={Guangqian Zhao},
pdfsubject={Affine ideal responses, source geometry, and covariance compensation}}
\theoremstyle{plain}
\newtheorem{maintheorem}{Theorem}

\newtheorem{theorem}{Theorem}[section]
\newtheorem{lemma}[theorem]{Lemma}
\newtheorem{proposition}[theorem]{Proposition}
\newtheorem{corollary}[theorem]{Corollary}
\theoremstyle{remark}
\newtheorem{example}[theorem]{Example}
\newtheorem*{remark}{Remark}
\numberwithin{equation}{section}
\newcommand{\E}{\mathbb E}
\newcommand{\Pbb}{\mathbb P}
\newcommand{\rank}{\operatorname{rank}}
\newcommand{\I}{\mathcal I}
\newcommand{\N}{\mathbb N}
\newcommand{\HS}{\mathrm{HS}}

\title[Deterministic principal-ideal responses]{Deterministic Principal-Ideal Responses\\
of Random Operators}
\author{Guangqian Zhao}
\address{School of Mathematical Sciences, University of Science and Technology of China, Hefei 230026, China}
\email{zhaoguangqian@mail.ustc.edu.cn}
\date{}
\begin{document}

\begin{abstract}
Let $X$ be a random compact operator for which each fixed compact shift
$X+B$ generates an almost surely deterministic principal ideal.
We show that all these ideals are determined by one ideal and one
coset. The ideal is generated almost surely by the difference of
independent copies of $X$, and the coset consists of all minimizing
shifts. For alternating quadratic evaluations under two-sided curvature
bounds, the spatial and kernel spectra determine the least ideal and
the covariance determines the coset. With the spatial weights and a
trace-class skew kernel fixed, the attainable raw ideals fill an explicit
order interval. Every ideal in the interval occurs with standard Gaussian
row marginals and joint covariance arbitrarily close to the identity.
For finite-dimensional operator spaces with infinite-rank coefficient
ideal we characterize almost-everywhere saturation in every translate
by a spectral summability criterion. For independent Gaussian rows we also
compute the least ideal attainable after observing a noise subspace.
For Wiener area, finitely many fixed Hilbert-valued temporal integrals
leave this ideal unchanged.
\end{abstract}

\maketitle

\section{Introduction}\label{intro:section}
Write $\I(T)$ for the two-sided ideal generated by a compact operator
$T$, with no norm closure, and $\I(u)=\I(\operatorname{diag}u)$.
Two operators generate the same principal ideal precisely when their
singular-value sequences dominate each other up to constants and finite
index dilations. Such an identity settles all two-sided ideal membership tests.

The response $B\mapsto\I(X+B)$ of a Borel random compact operator
$X$ is \emph{complete and deterministic} if, for each fixed compact
$B$, its value is almost surely a deterministic principal ideal.
Theorem~\ref{ideal:response-structure} proves that this holds exactly
when a deterministic principal ideal $J$ and compact operator $C$ satisfy
\begin{equation}\label{intro:affine-law}
 \I(X+B)=J+\I(B-C)\qquad\text{almost surely for each fixed }B.
\end{equation}
The pair $(J,C+J)$ is unique: independent copies satisfy
$\I(X-X')=J$ almost surely, and $C+J$ is the set of shifts attaining
the least ideal $J$. No moment assumption is needed.

For an alternating quadratic evaluation, the model data determine both
parts of the abstract response.
Let $Q\ge0$ be trace class with
$\rank Q\ge2$ and let $K\in\mathcal S_2$ be nonzero and skew-adjoint.
The centered array $Z$ has positive finite-coordinate densities
$ce^{-V}$ satisfying
\[
 \alpha I\le\nabla^2V\le\beta I,\qquad 0<\alpha\le\beta<\infty,
\]
uniformly over the coordinates and their number. The lower bound
controls moments; the upper bound gives conditional anti-concentration.
Put
\[
 q=s(Q),\quad a=\sqrt q,\quad\kappa=s(K),\qquad
 T_n=\sum_{j\ge n}\kappa_j^2,\quad
 E_n=\sum_{n\le i<j}q_iq_j,
\]
and $e_n^2=E_nT_n-E_{n+1}T_{n+1}$.
Let $\Gamma$ denote the full covariance of the array.
The centered evaluation $W_{Q,Z}(K)$ exists even when raw means
diverge. Theorem~\ref{ext:classification} gives
\begin{equation}\label{intro:model-law}
 \I(W_{Q,Z}(K)+B)=\I(e)+\I(B-\mathfrak C_\Gamma)
 \qquad\text{almost surely for each fixed compact }B.
\end{equation}
The exterior shell $e$ determines the least ideal.
The compensator $\mathfrak C_\Gamma$ is an absolutely convergent
sum of compressed covariance contractions, linear and continuous
in the full covariance for fixed spectral choices. Its class modulo
$\I(e)$ is intrinsic. Trace class is exactly the universal kernel
domain of raw evaluation in the curvature class.

For $K\in\mathcal S_1$, write $A_Z$ for the raw evaluation and put
$R_n=\sum_{j\ge n}q_j$, $S_n=\sum_{j\ge n}\kappa_j$ and
$p_n=R_nS_n-R_{n+1}S_{n+1}$. Keeping $Q,K$ fixed,
Theorem~\ref{ext:raw-interval} gives the entire range
\[
 \{L:\I(A_Z)=L\text{ almost surely for some curvature array }Z\}
 =\{\I(u):\I(e)\subseteq\I(u)\subseteq\I(p)\}.
\]
Every ideal in the interval occurs with standard Gaussian row marginals
and joint covariance arbitrarily close to the identity in operator norm.
The construction gives the minimizing shifts for each such law.

For independent standard Gaussian rows, write $A_Q(K)=W_{Q,Z}(K)$.
Theorem~\ref{ext:gaussian-information} determines the response to every
chosen compact correction measurable from observations on the range
of a fixed orthogonal projection $P$ on the noise space.
The one-sided residual kernel $(I-P)K$ determines the least attainable
ideal, and the negative conditional mean $-A_Q(PKP)$ determines the
coset of minimizing corrections. The projection may have infinite rank
and need not reduce $K$: the one-sided kernel retains interactions
between observed and unobserved coordinates.

For Wiener area this gives the covariance-tail response
\[
 \I(A_Q+D)=\I((R_n/n)_n)+\I(D)\qquad\text{almost surely}.
\]
It holds both for fixed compact $D$ and for compact corrections
measurable with respect to finitely many fixed Hilbert-valued temporal
integrals. The area ideal equals $\I(Q)$ exactly when
$R_n\le Cnq_n$. Examples distinguish the response from the unshifted
ideal and from all Schatten membership tests.

For the finite-dimensional existence problem, let $V$ be a real
operator space with basis
$A_1,\ldots,A_d$, and put
\[
 A(x)=\sum_{i=1}^d x_iA_i,\qquad
 S=\Bigl(\sum_{i=1}^d A_i^*A_i\Bigr)^{1/2},\qquad
 J=\I(S),\quad a_n=s_n(S).
\]
Assume that $J$ has infinite rank. The algebraic upper bound in
$B+V$ is $J+\I(B)$. Theorem~\ref{fs:main} proves that it is
attained for Lebesgue-almost every parameter in every translate
if and only if
\begin{equation}\label{intro:source-condition}
 \sum_{n\ge1}\frac{a_{mn}}{n a_n}<\infty
 \quad\text{for some integer }m\ge2.
\end{equation}
The exceptional sets then have Hausdorff dimension at most $d-1$,
sharply. If the condition fails, a dense $G_\delta$ set of compact
shifts in operator norm has no saturated parameter at all.

Condition~\eqref{intro:source-condition} depends only on $J$ and
expresses nuclearity of an intrinsic spectral-dilation scale.
It has an algebraic characterization: whenever
$J=\I(A)+\I(B)$, some nonzero $(s,t)\in\mathbb R^2$ satisfies
$\I(sA+tB)=J$ (Theorem~\ref{fs:intrinsic-dichotomy}).
Failure of nuclearity permits linear diagonal models inside $J$
that preserve subspace inclusions and sums in every finite dimension.

For pencils satisfying \eqref{intro:source-condition}, a source law gives
saturation for every pencil and shift exactly when it vanishes on all
sets in the class $\mathcal Z_d$ defined by affine-slab covers
(Theorem~\ref{fs:slab-source}). One trace-class ideal suffices to test this.

\smallskip\noindent\textbf{Related work.}
Principal-ideal calculus and finite dilations are classical
\cite{Calkin,KW,Simon}. Condition~\eqref{intro:source-condition}
is stronger than the soft-edged principal-ideal criterion
\cite[Proposition~4.6(ii)]{KW2007}; its scale interpretation uses
the classical nuclearity criterion \cite[Section~8]{Becnel}.
The polynomial extension uses Carbery--Wright
\cite[Theorem~8]{CW}. For $d\ge2$, $\mathcal Z_d$ is simultaneous
$(s,d-1)$-tube nullity for all $s\in(0,1]$
\cite[Section~4.5]{Fraser2026}; the contribution here is its operator
realization and source-law characterization. The exceptional-set theorem of \cite{LTV2022} concerns
mutual singularity of the singular parts of scalar spectral measures
under finite-rank self-adjoint perturbations.

The probability inputs include Brascamp--Lieb variance bounds
\cite{BL1976}, Harg\'e's convex comparison \cite{Harge2004},
and Gaussian comparison and concentration \cite{vH}.
The exterior-dilation estimates in \cite[Theorem~1]{ZhaoExterior}
give simultaneous Schatten bounds for $0<p<2$; our ideal identities
do not assert same-index equivalence of the full random spectrum.
Dickinson~\cite{Dickinson2007} studies optimal $L^2$ approximation
of iterated Brownian integrals from finitely many Gaussian linear
observations. For $Q$-Wiener processes,
Leonhard and R\"o\ss ler~\cite{LeonhardRossler2019} analyze simulation
algorithms given the increments and derive mean-square error bounds.
The observation theorem here determines the exact ideal response to
every chosen observation-measurable compact correction.

Throughout, probability-one events may depend on the fixed shift or
chosen correction. A prescribed countable family shares one event.

\section{Deterministic fixed-shift responses}\label{ideal:section}

All Hilbert spaces are real and separable. Write $\mathcal B,\mathcal K$ for
bounded and compact operators, with the operator norm. Singular values are
nonincreasing, include multiplicities and are zero-completed. They satisfy
\begin{equation}\label{var:approximation}
 s_n(T)=\inf_{\rank R<n}\|T-R\|.
\end{equation}
We compare rectangular operators after zero extension in a common
infinite-dimensional Hilbert space; this and complexification preserve nonzero
singular values. For $0<p<\infty$, write
$\mathcal S_p=\{T\in\mathcal K:\sum_ns_n(T)^p<\infty\}$ and
$\|T\|_p=(\sum_ns_n(T)^p)^{1/p}$; this is a quasi-norm when $p<1$.
In particular $\mathcal S_1$ and $\mathcal S_2$ are the trace-class and
Hilbert--Schmidt classes. We also write $\|T\|_{\HS}=\|T\|_2$ and
$(x\otimes y)z=x\langle y,z\rangle$. Unless indicated otherwise,
$\|T\|$ is the operator norm. We use \cite{Simon}
\begin{equation}\label{var:singularinequalities}
 s_n(ATB)\le\|A\|\|B\|s_n(T),\qquad
 s_{j+k-1}(S+T)\le s_j(S)+s_k(T).
\end{equation}
Thus $|s_n(S)-s_n(T)|\le\|S-T\|$. The ideal $\I(T)$ consists of all finite
sums $\sum_lA_lTB_l$ with bounded multipliers, and
$\I(a)=\I(D_a)$, where $D_a=\operatorname{diag}(a_n)$. No closure is taken. Nonzero finite-rank
operators all generate the finite-rank ideal.
Random compact operators are Borel random elements of the separable space
$\mathcal K$ with its operator norm; operator-valued expectations in
$L^p(\mathcal S_2)$ are Bochner expectations. Conditional laws of countable arrays are taken on their canonical standard
Borel coordinate spaces, so regular conditional distributions exist.
We use the standard conditional-expectation and martingale convergence
results of \cite[Chapters~8--9]{Kallenberg2021}.
A constant described as absolute is independent
of the dimensions and of the deterministic coefficients.

\begin{lemma}[Principal ideals and finite dilations]\label{var:idealcalculus}
For nonnegative nonincreasing null sequences $a,b$,
\begin{equation}\label{var:idealcriterion}
 \I(a)\subseteq\I(b)
 \quad\Longleftrightarrow\quad
 a_n\le Cb_{\lceil n/M\rceil}\quad(n\ge1)
 \quad\text{for some }C<\infty,\ M\in\mathbb N.
\end{equation}
Moreover, $\I(T)=\I(s(T))$ for every compact $T$, and
\begin{equation}\label{var:idealsum}
 \I(a)+\I(b)=\I(\max\{a,b\}),\qquad
 \I(a)\cap\I(b)=\I(\min\{a,b\}),
\end{equation}
where both extrema are coordinatewise. \end{lemma}
\begin{proof}
If $D_a=\sum_{j=1}^m A_jD_bB_j$, the singular-value inequalities give
\[
 s_{m(n-1)+1}(D_a)\le
 \sum_{j=1}^m\|A_j\|\|B_j\|b_n,
\]
which proves the forward implication in \eqref{var:idealcriterion}.
Conversely, for $(\Delta_m b)_n=b_{\lceil n/m\rceil}$, the isometries
$V_je_k=e_{m(k-1)+j}$ satisfy
$D_{\Delta_m b}=\sum_{j=1}^mV_jD_bV_j^*\in\I(b)$.
The bound $a\le C\Delta_m b$ then expresses $D_a$ as a bounded diagonal
multiple of $D_{\Delta_m b}$, with quotient zero at zero entries.
Singular-value decomposition gives $\I(T)=\I(s(T))$.

Splitting the coordinates according to $a_n\ge b_n$ proves the sum formula.
For the intersection, membership in both ideals bounds the singular
values by dilates of both sequences. Taking the larger dilation and
constant in \eqref{var:idealcriterion} gives a single bound by their minimum.
The converse follows from monotonicity.
These are classical consequences of Calkin's correspondence
\cite{Calkin}; see also \cite[Section~2]{KW}.
\end{proof}

\begin{lemma}[Completion from internal shifts]
\label{ideal:completion}\label{ideal:extraction}
Let $X$ be a random compact operator, $J$ a deterministic principal ideal,
and $C$ a deterministic compact operator. Suppose that $X+C\in J$
almost surely and that, for every fixed $D\in J$,
\[
 J\subseteq\I(X+C+D)\qquad\text{almost surely}.
\]
Then, for every fixed compact $B$,
\begin{equation}\label{ideal:completion-law}
 \I(X+B)=J+\I(B-C)\qquad\text{almost surely}.
\end{equation}
If $J$ has infinite rank, it suffices to test $D\in J$ with
$\|D\|<\varepsilon$, for any fixed $\varepsilon>0$.
\end{lemma}
\begin{proof}
Put $Y=X+C$. On the event $Y\in J$,
$\I(Y+D)+J=\I(D)+J$ for every compact $D$, by addition and subtraction.
It remains to prove $J\subseteq\I(Y+D)$ for each fixed $D$.
Write $J=\I(a)$ and choose a singular-value representation
$D=U\operatorname{diag}(b_n)V^*$, with $U,V$ partial isometries.
For integers $r\ge1$, put
\[
 D_r=U\operatorname{diag}
 \bigl(\min\{b_n,r\,a_{\lceil n/r\rceil}\}\bigr)V^*\in J.
\]
If $D\in J+\I(k)$, finite dilation gives some $r$ with
$b_n\le r(a_{\lceil n/r\rceil}+k_{\lceil n/r\rceil})$,
so $D-D_r\in\I(k)$ by coordinatewise domination.
Intersect the probability-one events for $Y\in J$ and these
countably many internal tests. For each sample set $K=\I(Y+D)$.
Since $D\in J+K$, some $r$ gives $Y+D_r\in K$, and hence
$J=\I(Y+D_r)\subseteq K$.

For the local version, write any $D\in J$ as a finite-rank operator
plus $D_0\in J$ with $\|D_0\|<\varepsilon$.
The local hypothesis gives $\I(Y+D_0)=J$; finite-rank perturbations
preserve the principal ideal of an infinite-rank compact operator.
\end{proof}
The lower bound may also be checked after a fixed bounded linear map
$\Lambda$ with $\I(\Lambda T)\subseteq\I(T)$: prove
$J\subseteq\I(\Lambda X+F)$ for every fixed compact $F$ and take
$F=\Lambda(C+D)$.

The infinite-rank restriction in the local version is necessary.
If $R\ne0$ has rank one and $Y=\pm R$ with equal probabilities,
then $\I(Y+D)=\mathcal F$ for every $D\in\mathcal F$ with
$\|D\|<\|R\|$, whereas $\I(Y-R)$ is zero or $\mathcal F$.
Here $\mathcal F$ is the ideal of finite-rank operators.

\begin{maintheorem}[Affine structure of deterministic fixed-shift responses]
\label{ideal:response-structure}
Let $X$ be a Borel random compact operator on a separable Hilbert space,
and let $X'$ be an independent copy. The following are equivalent:
\begin{enumerate}[label=\textup{(\roman*)}]
\item For every fixed compact $B$, $\I(X+B)$ is almost surely
      a deterministic principal ideal.
\item For every deterministic principal ideal $M$,
      $\mathbb P\{X-X'\in M\}\in\{0,1\}$; and for every fixed
      compact $B$, the ideals $\I(X+B),\I(X'+B)$ are almost surely
      comparable by inclusion.
\item There are a deterministic principal ideal $J$ and a deterministic
      compact operator $C$ such that, for every fixed compact $B$,
\begin{equation}\label{ideal:abstract-response}
 \I(X+B)=J+\I(B-C)\qquad\text{almost surely}.
\end{equation}
\item There is a deterministic compact $C$ such that, for every fixed
      compact $B$,
\begin{equation}\label{ideal:independent-absorption}
 \I(X+X'+C+B)\subseteq\I(X+B)\qquad\text{almost surely}.
\end{equation}
\end{enumerate}
When these conditions hold, $J$ and $C+J$ are unique,
\begin{equation}\label{ideal:abstract-difference}
 \I(X-X')=J\quad\text{almost surely},\qquad X+C\in J
 \quad\text{almost surely},
\end{equation}
and $J$ is the least response value, with minimizing shifts
\begin{equation}\label{ideal:minimizing-shifts}
 \{B:\I(X+B)=J\text{ almost surely}\}=C+J.
\end{equation}
\end{maintheorem}
\begin{proof}
Principal-ideal inclusion is a Borel relation on pairs of compact operators:
\[
 \{\I(S)\subseteq\I(T)\}
 =\bigcup_{c,m\in\N}\bigcap_{n\ge1}
   \{s_n(S)\le c\,s_{\lceil n/m\rceil}(T)\}.
\]
Equality is Borel as well. Thus two independent identically distributed
operators have the same principal ideal almost surely if and only if
that ideal is deterministic: for the nontrivial implication, fix one
full-measure section in Fubini's theorem. These arguments use no
measurable choice of ideal representatives.

Assume \textup{(ii)}. Every fixed coset of a principal ideal $M$
has probability zero or one: a coset of probability $p\in(0,1)$ would give
\[
 0<p^2\le\mathbb P\{X-X'\in M\}
 \le p^2+(1-p)^2<1,
\]
a contradiction. Fix $B$, let $\mu$ be the law of $X+B$, and set
\[
 q(y)=\mu\{z:\I(z)\subseteq\I(y)\},\qquad A=\{y:q(y)=0\}.
\]
The function $q$ is measurable, and the coset observation gives
$q(y)\in\{0,1\}$ for every $y$. On $A\times A$, inclusion in either
direction has product measure zero. Comparability therefore forces
$\mu(A)=0$. Hence two independent response ideals contain one another
almost surely, and the preceding Fubini argument proves \textup{(i)}.

Assume \textup{(i)} and take three independent copies $X_1,X_2,X_3$.
Conditioning on $X_1$ and changing signs gives
$\I(X_1-X_2)=\I(X_1-X_3)$ almost surely.
By Fubini, fix $x_0$ in place of $X_3$ so that
$\I(X_1-X_2)=\I(X_1-x_0)$ almost surely.
Condition \textup{(i)} at the fixed shift $-x_0$ makes the right side
a deterministic principal ideal $J$. Put $C=-x_0$.
Then $\I(X-X')=J$ and $\I(X+C)=J$ almost surely.
For any fixed $B$, let $L_B$ be its deterministic response. Since
\[
 J=\I(X_1-X_2)\subseteq\I(X_1+B)+\I(X_2+B)=L_B
 \qquad\text{almost surely},
\]
Lemma~\ref{ideal:completion} gives \textup{(iii)}.
The difference identity also proves \textup{(i)}$\Rightarrow$\textup{(ii)}.

Under \textup{(iii)}, $X'+C\in J\subseteq\I(X+B)$ almost surely,
which proves \textup{(iv)}. Conversely, put $Y=X+C$ and
$Y'=X'+C$. Condition \textup{(iv)} at $B=C$ gives
$Y'\in\I(Y)$ almost surely. Exchangeability gives the reverse
inclusion, so $\I(Y)=J$ almost surely for a deterministic principal
ideal $J$. For any fixed $B$, subtracting $X+B$ in
\eqref{ideal:independent-absorption} yields
$Y'\in\I(X+B)$ and therefore $J\subseteq\I(X+B)$ almost surely.
Lemma~\ref{ideal:completion} again proves \textup{(iii)}.
The implication \textup{(iii)}$\Rightarrow$\textup{(i)} is immediate.

The value in \eqref{ideal:abstract-response} is always at least $J$
and equals $J$ exactly when $B-C\in J$. This proves uniqueness and
\eqref{ideal:minimizing-shifts}. The compensations in \textup{(iv)}
are precisely this same coset.
\end{proof}

Neither part of \textup{(ii)} can be dropped. The rank-one example
above has comparable but random response ideals. Conversely, the Gaussian
plane $X$ in Theorem~\ref{fs:intrinsic-dichotomy} satisfies
$X-X'\overset{\mathrm d}{=}\sqrt2X$.
Proposition~\ref{fs:response-defects} therefore gives its zero-one
difference tests, but independent unshifted response ideals are almost
surely incomparable.

For any deterministic compact $m$ and algebraic ideal $L$,
\eqref{ideal:abstract-response} at shift $-m$ gives
$X-m\in L$ almost surely exactly when $J+\I(m+C)\subseteq L$.
Containment is understood on a measurable probability-one event;
no Borel assumption on $L$ is needed.
Thus $-C+J$ is the least deterministic affine ideal carrier, and,
for $J\subseteq L\subseteq\mathcal K$,
\begin{equation}\label{ideal:affine-ball}
 \{B\in\mathcal K:\I(X+B)\subseteq L\text{ almost surely}\}=C+L.
\end{equation}

For an independent random compact $Y$, the response formula and Fubini give
$\I(X+Y)=J+\I(Y-C)$ almost surely. A block sum of compact generators
represents the ideal sum, so the relevant comparison is Borel.
Consequently, for finitely many independent $X_i$ with pairs
$(J_i,C_i+J_i)$, nonzero real $t_i$, and each fixed compact $B$,
\begin{equation}\label{ideal:independent-linear-combination}
 \I\left(\sum_i t_iX_i+B\right)
 =\sum_iJ_i+\I\left(B-\sum_i t_iC_i\right)
 \qquad\text{almost surely}.
\end{equation}
If $X$ has the same law as $-X$, uniqueness gives $C+J=-C+J$,
so $C\in J$ and $\I(X+B)=J+\I(B)$, with $J=\I(X)$ almost surely.
These conclusions assume complete determinacy; the next section
gives its finite-source criterion.

\section{A geometric dichotomy for finite operator spaces}
\label{source:section}\label{fs:section}

Let $V$ be a $d$-dimensional real subspace of $\mathcal K(H)$,
$d\ge1$, with basis $A_1,\ldots,A_d$, and put
\begin{equation}\label{fs:column-profile}
 \begin{gathered}
 A(x)=\sum_{i=1}^d x_iA_i,\qquad
 S=\Bigl(\sum_{i=1}^d A_i^*A_i\Bigr)^{1/2},\\
 a_n=s_n(S),\qquad J=\sum_{i=1}^d\I(A_i)=\I(S).
 \end{gathered}
\end{equation}
Indeed, the compact column operator $\mathcal A:H\to H^d$,
$\mathcal A u=(A_1u,\ldots,A_du)$, has polar decomposition
$\mathcal A=US$. Writing $U_i$ for the components of the contraction
$U$ gives $A_i=U_iS$ and $S=U^*\mathcal A=\sum_iU_i^*A_i$,
which prove the ideal equality and $s_n(A_i)\le a_n$.
An invertible real change of coefficient basis multiplies
$\mathcal A$ by a bounded invertible matrix on $H^d$; hence it changes
$a_n$ by at most fixed multiplicative constants at the same index.
Assume that $J$ contains an infinite-rank operator, so that
$a_n>0$ and $a_n\downarrow0$. For a fixed compact perturbation $B$, define
\[
 E_B=\{x\in\mathbb R^d:\I(A(x)+B)\ne J+\I(B)\}.
\]
Thus $E_B$ records failure to attain the algebraic upper bound
$J+\I(B)$ on the affine space $B+V$. Its dimension does not depend on the
choice of basis of $V$.

\begin{maintheorem}[Geometric saturation dichotomy]\label{fs:main}
The following conditions are equivalent:
\begin{enumerate}[label=\textup{(\roman*)}]
\item For some integer $m\ge2$,
\begin{equation}\label{fs:criterion}
 \sum_{n\ge1}\frac{a_{mn}}{n a_n}<\infty.
\end{equation}
\item For every compact $B$, $\dim_{\mathrm H}E_B\le d-1$.
\item For every compact $B$, $\I(A(x)+B)=J+\I(B)$ for
Lebesgue-almost every $x\in\mathbb R^d$.
\item For every compact $B$, the set $\mathbb R^d\setminus E_B$
is dense in $\mathbb R^d$.
\item For independent standard Gaussian $G,G'$ in $\mathbb R^d$
and every fixed compact $B$, the ideals $\I(A(G)+B)$ and
$\I(A(G')+B)$ are almost surely comparable by inclusion.
\end{enumerate}
When these conditions hold, every $E_B$ is a $G_\delta$ set and
$\mathcal H^q(E_B)=0$ for every $q>d-1$; the dimension bound is sharp.
When they fail, there is a dense $G_\delta$ set of compact $B$
in operator norm such that
\[
 J\not\subseteq\I(A(x)+B)\qquad\text{for every }x\in\mathbb R^d.
\]
For each such $B$, every principal-ideal response fiber lies in a
proper affine hyperplane. Under any source law that gives hyperplanes
zero mass, two independent responses are almost surely incomparable.
Moreover, for every $R,\varepsilon>0$
there is $D\in J$ with $\|D\|<\varepsilon$ for which
$\I(A(x)+D)\subsetneq J$ throughout $[-R,R]^d$.
\end{maintheorem}

The spectral condition depends only on $J$. No commutation or
spectral comparability is assumed, and nonzero directions in $V$
may have finite rank. The residual alternative is a statement in
the norm topology of perturbations, not a probability assertion
about them. Below, $d$ denotes the number of displayed coefficients
and may exceed the dimension of their span. Linear independence is imposed
only where stated.

\subsection{Finite sampling and the order of response ideals}

\begin{proposition}[Finite-sample envelopes and diffuse response defects]
\label{fs:response-defects}
Let $X=A_0+A(\xi)$, where $A(x)=\sum_{i=1}^d x_iA_i$ and
the law of $\xi$ gives zero mass to every proper affine hyperplane.
Put $J=\sum_i\I(A_i)$ and $F_B=J+\I(A_0+B)$.
For every fixed algebraic two-sided ideal $M$ and compact $B$,
\[
 \mathbb P\{X+B\in M\}=\mathbf1_{\{F_B\subseteq M\}}.
\]
For $d+1$ independent copies $X_0,\ldots,X_d$, on one
probability-one event,
\[
 \sum_{j=0}^d\I(X_j+B)=F_B
 \qquad\text{for every compact }B.
\]
For each fixed $B$, set $L_B=\I(X+B)$ and
$p_B=\mathbb P\{L_B=F_B\}$, and let $L'_B$ be an independent copy.
The only principal ideal that can have positive probability as $L_B$
is $F_B$, and
\[
 \mathbb P\{L_B\subseteq L'_B\}=p_B,
 \qquad \mathbb P\{L_B=L'_B\}=p_B^2.
\]
Thus two deficient independent responses are almost surely incomparable,
and $L_B$ is deterministic if and only if $L_B,L'_B$ are almost surely
comparable.
\end{proposition}
\begin{proof}
For fixed $M,B$, the set
$\{x:A_0+A(x)+B\in M\}$ is empty or an affine subspace.
It is all of $\mathbb R^d$ exactly when $F_B\subseteq M$;
otherwise it lies in a proper affine hyperplane. This proves the first
assertion, including measurability even for an arbitrary algebraic $M$.
The source vectors $\xi_0,\ldots,\xi_d$ are almost surely affinely
independent. On that event their differences recover every $A_i$,
and their shifted values then recover $A_0+B$, for all $B$ at once.
This proves the finite-sample identity.

Every $L_B$ is contained in $F_B$. Conditional on the second sample,
the first assertion applied to $M=L'_B$ shows that the conditional
probability of $L_B\subseteq L'_B$ is
$\mathbf1_{\{L'_B=F_B\}}$. The ideal-comparison relation is Borel,
so Fubini applies. Equality therefore occurs, up to a null event,
exactly when both responses equal $F_B$. Independence gives the two
formulas. A fixed positive-probability ideal value must be $F_B$ by
the first assertion. Finally the probability of comparability is
$2p_B-p_B^2$, which equals one exactly when $p_B=1$.
\end{proof}

\subsection{Spectral dilation and a deterministic slab bound}

By Lemma~\ref{var:idealcalculus}, inclusion of principal ideals is
equivalent to a singular-value comparison with a finite index dilation.
For an infinite-rank target $T$, equivalently one may use
$s_{Mn}(S)\le C s_n(T)$ after changing these constants.
These countable singular-value comparisons also show that every
fixed principal-ideal response fiber is Borel.

\begin{lemma}[Dilation normalization]\label{fs:normalization}
For a positive decreasing null sequence $a$, \eqref{fs:criterion}
is equivalent to
\begin{equation}\label{fs:dyadic}
 \exists L\in\mathbb N:\qquad
 \sum_{j\ge0}\frac{a_{2^{j+L}}}{a_{2^j}}<\infty.
\end{equation}
It is invariant under equality of principal ideals. When it holds,
for every $\eta>0$ there is $K\ge1$ such that
$\sum_j(a_{2^{j+K}}/a_{2^j})^\eta<\infty$.
\end{lemma}
\begin{proof}
Monotonicity gives
\[
 \frac12\frac{a_{m2^{j+1}}}{a_{2^j}}
 \le\sum_{n=2^j}^{2^{j+1}-1}\frac{a_{mn}}{n a_n}
 \le\frac{a_{m2^j}}{a_{2^{j+1}}}.
\]
Choose $2^L\ge2m$ for one implication and $m=2^{L+1}$ for the other.
If $\I(a)=\I(b)$, put $\alpha_j=a_{2^j}$, $\beta_j=b_{2^j}$.
There are integers $P,Q$ and constants $C,D$ such that
$\alpha_{j+P}\le C\beta_j$ and $\beta_{j+Q}\le D\alpha_j$.
Hence
\[
 \frac{\beta_{j+P+Q+L}}{\beta_j}
 \le CD\frac{\alpha_{j+P+L}}{\alpha_{j+P}},
\]
which transfers summability; reverse the roles for equivalence.
Finally set $r_L(j)=\alpha_{j+L}/\alpha_j$. Since
$r_{kL}(j)=\prod_{h=0}^{k-1}r_L(j+hL)$, H\"older gives
\[
 \sum_jr_{kL}(j)^{1/k}
 \le\prod_{h=0}^{k-1}
       \left(\sum_jr_L(j+hL)\right)^{1/k}<\infty.
\]
Choose $1/k\le\eta$ and use $0<r_{kL}(j)\le1$.
\end{proof}

A slab of width $w$ means a set between parallel hyperplanes
at distance $w$ in $\mathbb R^d$.

\begin{lemma}[Thin slabs of bad parameters]\label{fs:slab}
For every compact $B$, $n\ge1$ and $t>0$, the set
\[
 F_{n,t}(B)=\{x:s_n(A(x)+B)\le t\}
\]
is contained in a slab of width at most $4dt/a_{2dn}$.
\end{lemma}
\begin{proof}
First intersect the set with a closed bounded cube, obtaining a compact
set $F$. If its affine span is proper, a hyperplane contains it.
Otherwise choose a maximum-volume simplex $x_0,\ldots,x_d$ in $F$
and set $M=[x_1-x_0\ \cdots\ x_d-x_0]$.
Replacing any vertex by a point of $F$ shows that all its barycentric
coordinates have absolute value at most one. Consequently
$F\subseteq x_0+M[-1,1]^d$.
Put $Y_j=A(x_j)+B$, and choose $L_j$ of rank at most $n-1$
with $\|Y_j-L_j\|\le t$ for $0\le j\le d$.
The column operator $\mathcal V=(Y_j-Y_0)_{j=1}^d$ has the
approximant $(L_j-L_0)_{j=1}^d$ of rank at most $2d(n-1)$,
with error at most $2\sqrt d\,t$.
Since $\mathcal A=(M^{-\mathsf T}\otimes I)\mathcal V$ and
$s_k(\mathcal A)=s_k(S)=a_k$, it follows that
\[
 a_{2dn}\le a_{2d(n-1)+1}
       \le\frac{2\sqrt d\,t}{\sigma_{\min}(M)}.
\]
If $u$ is a left singular unit vector for $\sigma_{\min}(M)$,
the width of $M[-1,1]^d$ in direction $u$ is at most
\[
 2\|M^*u\|_1\le2\sqrt d\,\sigma_{\min}(M)
                    \le4dt/a_{2dn}.
\]
This proves the bound on every bounded cube.
If the full set is nonempty, fix a point in it and increase the cubes.
The slab normals have a convergent subsequence; their offsets remain
bounded because each slab contains the fixed point and has bounded width.
A limiting slab contains the full set.
\end{proof}

\begin{lemma}[Slab covers of exceptional parameters]
\label{fs:hausdorff}
For any finite coefficient list as above satisfying
\eqref{fs:criterion}, $E_B$ is a $G_\delta$ set.
For every $\eta>0$ and every fixed compact $B$, there are open
slabs $S_j$ such that
\[
 E_B\subseteq\limsup_jS_j,\qquad
 \sum_j\operatorname{width}(S_j)^\eta<\infty.
\]
Consequently $\mathcal H^q(E_B)=0$ for every $q>d-1$.
The dimension bound is sharp even for a linearly independent list
whose nonzero linear combinations all have infinite rank.
\end{lemma}
\begin{proof}
Since $A(x)\in J$, subtraction shows that failure of the defining
response identity is equivalent to $J\not\subseteq\I(A(x)+B)$.
The finite-dilation criterion therefore gives
\[
 E_B=\bigcap_{c,m\in\mathbb N}\ \bigcup_{n\ge1}
 \{x:a_n>c\,s_{\lceil n/m\rceil}(A(x)+B)\},
\]
which is $G_\delta$ by continuity of singular values.

Fix $\eta>0$ and an integer $C$ with $2^C\ge2d$.
By Lemma~\ref{fs:normalization}, choose $K$ such that
\[
 \sum_{j\ge0}
 \left(\frac{a_{2^{j+C+K}}}{a_{2^{j+C}}}\right)^\eta<\infty.
\]
Put $F_j=\{x:s_{2^j}(A(x)+B)\le a_{2^{j+C+K}}\}$.
Lemma~\ref{fs:slab}, followed by a slight enlargement, places $F_j$
in an open slab $S_j$ of width at most
\[
 w_j=8d\,\frac{a_{2^{j+C+K}}}{a_{2^{j+C}}}.
\]
Outside $\limsup_jF_j$, dyadic interpolation gives
\[
 s_n(A(x)+B)\ge a_{2^{C+K+1}n}
 \quad\text{for all sufficiently large }n.
\]
The output has infinite rank, so its finite head can be absorbed into
a constant. Finite dilation and subtraction now give
$\I(A(x)+B)=J+\I(B)$.
Thus $E_B\subseteq\limsup_jF_j\subseteq\limsup_jS_j$,
proving the slab-cover assertion.

For $q>d-1$, apply this assertion with $\eta=q-d+1$.
Inside a fixed bounded cube, a slab of width $w_j\le1$ is covered
by $O(w_j^{-(d-1)})$ balls of radius $O(w_j)$.
Their $q$th diameter powers have total $O(w_j^\eta)$.
Summability and then exhaustion by cubes give $\mathcal H^q(E_B)=0$.

For sharpness, take orthogonal infinite-dimensional blocks and let
$A_1$ have singular values $e^{-n}$ on the first block, while each
$A_i$, $i\ge2$, has singular values $e^{-n^2}$ on block $i$.
Then $J=\I((e^{-n})_n)$ satisfies \eqref{fs:criterion}, and
$E_0=\{x_1=0\}$. Every nonzero source direction has infinite rank.
\end{proof}

\begin{remark}[Finite-rank coefficients]
If $J$ is nonzero and has finite rank, then $E_B$ is empty for
infinite-rank $B$, and
\[
 E_B=\{x:A(x)+B=0\}
\]
for finite-rank $B$; the latter set is empty or a proper affine
subspace. Indeed, finite-rank perturbations preserve the principal
ideal of an infinite-rank operator. If $J=0$, every $E_B$ is empty.
\end{remark}

\subsection{Local and residual failure}

The converse uses a deterministic block construction valid for arbitrary
noncommuting coefficients. We record its two elementary selection steps.

\begin{lemma}[Power-divergent disjoint blocks]\label{fs:blocks}
Let $b_j>0$ decrease to zero and suppose
$\sum_j b_{j+L}/b_j=\infty$ for every integer $L\ge1$.
Then $\sum_j(b_{j+L}/b_j)^p=\infty$ for every $p>0$ and $L\ge1$.
Given $p>0$, $L\ge1$, $0<\delta<1$, a cutoff, and $R>0$, there are finitely
many ordered disjoint integer intervals $[u_i,v_i)$ beyond that cutoff,
with $v_i>u_i+L$, such that
\begin{equation}\label{fs:block-budget}
 b_{v_i}\le\delta b_{u_i+L},\qquad
 \sum_i\left(\delta\frac{b_{u_i+L}}{b_{u_i}}\right)^p\ge R.
\end{equation}
\end{lemma}
\begin{proof}
If the power series converged for some $p,L$, take an integer $k\ge p$.
The product identity in Lemma~\ref{fs:normalization} and H\"older give
$\sum_j(b_{j+kL}/b_j)^{p/k}<\infty$.
Since all ratios lie in $(0,1]$, this contradicts divergence at $kL$.

For the selection, one residue class modulo $L$ has divergent power sum.
Write $c_k=b_{r_0+kL}$ and $t_k=c_{k+1}/c_k$ on this class.
If $\limsup t_k>0$, choose infinitely many $t_k\ge\zeta>0$.
After choosing $k$, take $\ell>k+1$ with
$c_\ell\le\delta c_{k+1}$, retain $[k,\ell)$, and choose the
next good index at least $\ell$. The resulting disjoint intervals
each have weight at least $(\delta\zeta)^p$.
If $t_k\to0$, then eventually $t_{k+1}\le\delta$.
One parity of $\sum_k t_k^p$ diverges; its intervals $[k,k+2)$
are disjoint and have the required endpoint drop.
In either case a finite selection reaches total weight $R$.
Returning to the original indices proves \eqref{fs:block-budget}.
\end{proof}

\begin{lemma}[Covering by prescribed radii]\label{fs:cover}
If finitely many $0<\rho_i\le1/2$ satisfy $\sum_i\rho_i^d\ge1$,
there are $z_i\in[0,1]^d$ whose closed sup-norm balls of radii
$\rho_i$ cover $[0,1]^d$.
\end{lemma}
\begin{proof}
Let $\ell_i$ be the largest dyadic number at most $2\rho_i$.
Then $\rho_i<\ell_i\le2\rho_i$ and $\sum_i\ell_i^d\ge1$.
Process these side lengths in decreasing order, placing dyadic cubes
inside $[0,1]^d$ with disjoint interiors until the unit cube is full.
At each step the unfilled region is a union of grid cubes of the
current side length. Unless already full, a placement is possible.
The total prescribed volume forces completion before the list ends.
Each placed cube lies in the corresponding sup-norm ball about its
center; give arbitrary centers to unused radii.
\end{proof}

\begin{lemma}[Arbitrarily small perturbations losing a whole cube]
\label{fs:failure}
Suppose \eqref{fs:criterion} fails. For every $R,\varepsilon>0$
there is $B\in J$, $\|B\|<\varepsilon$, such that
$\I(A(x)+B)\subsetneq J$ for every $x\in[-R,R]^d$.
Every response fiber restricted to this cube lies in an affine
hyperplane. Consequently, for any probability law $\mu$ assigning
zero mass to affine hyperplanes, an arbitrarily small $B\in J$
can make $\I(A(\xi)+B)$, $\xi\sim\mu$, nondeterministic.
If $A_i^*=\epsilon A_i$ for a common $\epsilon\in\{1,-1\}$,
$B$ can be chosen with $B^*=\epsilon B$.
\end{lemma}
\begin{proof}
Homogeneity reduces the construction to $R=1$:
$A(x)+RB=R(A(x/R)+B)$, and the unit-cube perturbation may be chosen
with norm below $\varepsilon/R$.
Use $S$ and $a_n=s_n(S)$ from \eqref{fs:column-profile}.
Let $F_n$ project onto the first $n-1$ eigenvectors of $S$
corresponding to its positive eigenvalues. Then
\begin{equation}\label{fs:common-projections}
 \rank F_n=n-1,\qquad
 \|A_i(I-F_n)\|\le\|S(I-F_n)\|=a_n.
\end{equation}
The projections are nested; no commutation between $A_i$ and $S$
is required.

At stage $h\ge1$, set $\delta_h=2^{-h}/(2d)$.
Apply Lemma~\ref{fs:blocks} to $b_j=a_{2^j}$ with $p=d$, $L=h$,
and total budget $1$, placing its blocks beyond all earlier stages.
Put $n_\ell=2^{u_\ell}$, $N_\ell=2^{v_\ell}$ and
$\rho_\ell=\delta_h a_{2^h n_\ell}/a_{n_\ell}$.
At each stage,
\[
 a_{N_\ell}\le\delta_h a_{2^h n_\ell},\qquad
 \sum_\ell\rho_\ell^d\ge1,\qquad0<\rho_\ell<1/2.
\]
Lemma~\ref{fs:cover}, followed by rescaling, gives centers
$z_\ell\in[-1,1]^d$ whose sup-norm balls of radii $2\rho_\ell$
cover that cube. The intervals $[n_\ell,N_\ell)$ are ordered
and disjoint across all stages. Define
\[
 P_\ell=F_{N_\ell}-F_{n_\ell},\qquad
 D_i=\sum_\ell z_{\ell,i}P_\ell,
 \qquad B=-\sum_{i=1}^d A_iD_i.
\]
The $P_\ell$ are pairwise orthogonal projections, so each series
converges strongly to a self-adjoint $D_i$ with $\|D_i\|\le1$.
The finite sum defining $B$ is compact and belongs to the algebraic
ideal $J$. Choose the initial cutoff $n_0$ with
$d a_{n_0}<\varepsilon$. Every $D_i$ vanishes on $F_{n_0}$,
and \eqref{fs:common-projections} gives $\|B\|<\varepsilon$.

Fix $x\in[-1,1]^d$ and at each stage choose a block $\ell$ with
$\|x-z_\ell\|_\infty\le2\rho_\ell$.
Since $D_i$ commutes with $S$, its multiplier $x_iI-D_i$ has
norm at most $2\rho_\ell$ on $P_\ell H$ and at most $2$ afterwards.
Using $A_i=U_iS$ and $\|U_i\|\le1$, we obtain
\[
 \|(A(x)+B)(I-F_{n_\ell})\|
 \le2d\max\{\rho_\ell a_{n_\ell},a_{N_\ell}\}
 \le2^{-h}a_{2^h n_\ell}.
\]
Deleting rank at most $n_\ell-1$ yields
\begin{equation}\label{fs:bad-spectrum}
 s_{n_\ell}(A(x)+B)\le2^{-h}a_{2^h n_\ell}.
\end{equation}
If $J\subseteq\I(A(x)+B)$, finite dilation would give
$a_{Mk}\le C s_k(A(x)+B)$ with fixed $M,C$ for this $x$.
Choose $h$ so large that $2^h\ge M$ and $C2^{-h}<1$.
At $k=n_\ell$, monotonicity gives
$a_{Mk}\ge a_{2^h n_\ell}$, contradicting \eqref{fs:bad-spectrum}.
Since $A(x)+B\in J$, every response on the cube is a proper subideal of $J$.

For the symmetry assertion, replace $B$ by
$B_\epsilon=(B+\epsilon B^*)/2$.
This belongs to $J$, has norm at most $\|B\|$, and satisfies
$B_\epsilon^*=\epsilon B_\epsilon$.
For each selected block, put $Q=I-F_{n_\ell}$ and
$Y=A(x)+B$. Since
$A(x)+B_\epsilon=(Y+\epsilon Y^*)/2$, the ordinary tail estimate gives
\[
 \|Q(A(x)+B_\epsilon)Q\|
 \le\|QYQ\|\le\|YQ\|\le2^{-h}a_{2^h n_\ell}.
\]
Double-sided deletion costs rank at most $2(n_\ell-1)$.
Thus the preceding contradiction applies at
$k=2n_\ell-1\le2n_\ell$, taking $2^h\ge2M$.

For either resulting shift $\widetilde B$, a fiber with ideal value
$K\subsetneq J$ lies in $\{x:A(x)+\widetilde B\in K\}$, a proper
affine subspace by the observation in Proposition~\ref{fs:response-defects}.
For the probabilistic assertion, choose $R$ with $\mu([-R,R]^d)>0$
and use the construction on that cube. An almost surely constant
response would give a positive-mass restricted fiber, contradicting
hyperplane-nullity. Hence the response is nondeterministic.
\end{proof}

\begin{proposition}[Residual loss for every parameter]
\label{fs:generic-failure}
Let $A(x)=\sum_{i=1}^d x_iA_i$, let
$J=\sum_i\I(A_i)=\I(a)$ have infinite rank, and suppose that
\eqref{fs:criterion} fails. There is a dense $G_\delta$ set
$\mathcal G\subseteq\mathcal K(H)$ such that
\begin{equation}\label{fs:global-loss}
 J\not\subseteq\I(A(x)+B)
       \qquad(B\in\mathcal G,\ x\in\mathbb R^d).
\end{equation}
Every principal-ideal response fiber
$\{x:\I(A(x)+B)=K_0\}$, $B\in\mathcal G$, lies in an affine
hyperplane. Consequently, for every probability law $\mu$ assigning
zero mass to affine hyperplanes, every such fiber has $\mu$-measure zero.

If all coefficients satisfy $A_i^*=\epsilon A_i$ for one
$\epsilon\in\{1,-1\}$, the same assertions hold with a dense
$G_\delta$ set in the closed real Banach space
$\mathcal K_\epsilon(H)=\{B\in\mathcal K(H):B^*=\epsilon B\}$.
\end{proposition}
\begin{proof}
Work in $\mathcal K(H)$, or in $\mathcal K_\epsilon(H)$ for the
symmetry assertion. For an integer $R\ge1$, put
\[
 \mathcal E_R=\{B:J\not\subseteq\I(A(x)+B)
                       \text{ for every }x\in[-R,R]^d\}.
\]
This set is dense. Indeed, given $B_0$ and $\eta>0$, choose
finite-rank $F$ with $\|B_0-F\|<\eta/2$, preserving symmetry
when required. Lemma~\ref{fs:failure} supplies $D\in J$ with
$\|D\|<\eta/2$ such that $\I(A(x)+D)\subsetneq J$ throughout
the cube, again with the required symmetry.
Writing $\mathcal F$ for the finite-rank ideal, for each
$T=A(x)+D$ we have
\[
 \I(T+F)\subseteq\I(T)+\mathcal F\subsetneq J.
\]
The last inclusion holds because every nonzero operator ideal
contains $\mathcal F$, while $\mathcal F\subsetneq J$ also handles
$T=0$. Thus $F+D\in\mathcal E_R$ and $\|F+D-B_0\|<\eta$.

We next verify that $\mathcal E_R$ is $G_\delta$. For positive
integers $M,C,N$, put
\begin{align*}
 f_{R,M,C,N}(B)&=
 \max_{x\in[-R,R]^d}\min_{1\le n\le N}
                  \{C s_n(A(x)+B)-a_{Mn}\},\\
 \mathcal O_{R,M,C}&=\bigcup_{N\ge1}\{B:f_{R,M,C,N}(B)<0\}.
\end{align*}
Each $f_{R,M,C,N}$ is $C$-Lipschitz in operator norm, so
$\mathcal O_{R,M,C}$ is open in the chosen space. Finite dilation
gives
\[
 \mathcal E_R=\bigcap_{M,C\ge1}\mathcal O_{R,M,C}.
\]
For clarity, the pointwise condition on the left says that for
every $M,C$ and every $x$ some index $n$ gives a strict negative
value. Continuity supplies an open cover of the compact cube by
these witnessing inequalities. A finite subcover supplies one $N$;
the continuous minimum then has a strictly negative maximum.
This proves the displayed equality with its stated quantifiers.

Baire's theorem makes $\mathcal G=\bigcap_{R\ge1}\mathcal E_R$
a dense $G_\delta$ set, and its defining property is
\eqref{fs:global-loss}. If one response fiber contained $d+1$
affinely independent points, subtracting their operator values and
inverting the coefficient matrix would put every $A_i$ in that
response ideal, contradicting \eqref{fs:global-loss}.
Thus each fiber lies in an affine hyperplane, as asserted.
\end{proof}

\begin{proof}[Proof of Theorem~\ref{fs:main}]
Lemma~\ref{fs:hausdorff} proves \textup{(i)}$\Rightarrow$\textup{(ii)}
and the sharper Hausdorff-measure assertion. The implications
\textup{(ii)}$\Rightarrow$\textup{(iii)}$\Rightarrow$\textup{(iv)}
follow because $d$-dimensional Hausdorff-null sets are Lebesgue-null,
and every nonempty open set has positive Lebesgue measure.
If \textup{(i)} fails, Lemma~\ref{fs:failure} gives a fixed
$B\in J$ for which an entire cube misses the saturated response,
contradicting \textup{(iv)}. Proposition~\ref{fs:response-defects}
identifies \textup{(v)} with almost-sure saturation for standard Gaussian
parameters, equivalently \textup{(iii)} since Gaussian and Lebesgue
measure have the same null sets. The remaining assertions follow from
Lemma~\ref{fs:failure}, Proposition~\ref{fs:generic-failure}, and
Proposition~\ref{fs:response-defects}.
\end{proof}

For self-adjoint or skew-adjoint coefficient spaces, the same
geometric equivalence holds when perturbations are restricted to
the corresponding closed real subspace: sufficiency is unchanged,
and the relative residual obstruction proves necessity.
One cannot generally strengthen global loss by requiring $B\in J$.
For $T=\operatorname{diag}(n^{-p})$ and $B\in\I(T)$,
finite dilation gives $s_n(B)\le C_Bn^{-p}$, whence
\[
 s_n(tT+B)\ge |t|s_{2n-1}(T)-s_n(B)
 \ge (|t|2^{-p}-C_B)n^{-p}.
\]
Thus $\I(tT+B)=\I(T)$ for all sufficiently large $|t|$.

\subsection{Nuclearity and scalar generation}

A spectral dilation by $2^p$ becomes a shift by $p$ on dyadic indices.
Put $\alpha_j=a_{2^j}$ for $j\ge0$ and, for integers $p\ge0$, define
\begin{equation}\label{fs:scale}
 \mathscr H_p(a)=\left\{z:
    \|z\|_{a,p}^2=\sum_{j\ge0}\frac{|z_j|^2}{\alpha_{j+p}^2}<\infty
                    \right\},\qquad
 \mathscr E(J)=\bigcap_{p\ge0}\mathscr H_p(a).
\end{equation}
The increasing Hilbert norms give a complete metrizable locally convex
topology; finitely supported sequences are dense in each Hilbert completion.

\begin{proposition}[Nuclear interpretation]\label{fs:nuclear}
The space $\mathscr E(J)$ is independent of the decreasing positive
generator $a$ up to equivalence of its defining norms, with the identity
on sequence coordinates. It is nuclear if and only if
\eqref{fs:criterion} holds.
\end{proposition}
\begin{proof}
For another generator $b$, Lemma~\ref{fs:normalization} gives
$\|z\|_{a,p}\le D\|z\|_{b,p+Q}$ and
$\|z\|_{b,p}\le C\|z\|_{a,p+P}$.
Thus the two spaces coincide and have equivalent defining norms.
The inclusion $\mathscr H_{p+L}(a)\to\mathscr H_p(a)$ has diagonal
entries $\alpha_{j+p+L}/\alpha_{j+p}$ under the natural isometries
to $\ell^2$. Under \eqref{fs:dyadic}, the same $L$ makes every such
inclusion trace class: its entry sum is the tail starting at $p$
of the summable dyadic series.
Conversely, the countably-Hilbert nuclearity criterion
\cite[Definition~8.1 and Remark~8.2]{Becnel} supplies a trace-class
inclusion from a stronger Hilbert norm into $\mathscr H_0(a)$,
which is precisely \eqref{fs:dyadic}.
If nuclearity is formulated using Hilbert--Schmidt inclusions,
compose two such inclusions to obtain trace class.
\end{proof}

Nuclearity of $\mathscr E(J)$ is a condition on the connecting maps
of its Hilbert scale; it differs from trace-class membership of a generator. For comparison, the classical soft-edged
principal-ideal criterion is $a_{mn}/a_n\to0$ for some $m$
\cite[Proposition~4.6(ii)]{KW2007}; it corresponds here to a compact
inclusion between two such Hilbert spaces. The response criterion
requires a trace-class inclusion.

For example, $a_n=n^{-p}$, $p>0$, fails \eqref{fs:criterion}, whereas
$a_n=e^{-n^\delta}$ satisfies it for every $\delta>0$.
There is also a strict gap within rapidly decreasing spectra.
For a sufficiently large fixed $c$ and $\gamma>0$, set
\[
 a_n=\exp\{-\log(n+c)[\log\log(n+c)]^\gamma\}.
\]
These sequences belong to every $\ell^p$, $p>0$, and are soft-edged,
but \eqref{fs:criterion} holds exactly when $\gamma\ge1$.
Indeed a fixed dyadic ratio has logarithm
$-L\log2\,(\log j)^\gamma(1+o(1))$.
For $\gamma=1$ the ratio is asymptotic to a positive constant times
$j^{-L\log2}$; for $0<\gamma<1$ every fixed-dilation ratio tends
to zero but is nonsummable.

\begin{theorem}[An intrinsic two-generator dichotomy]
\label{fs:intrinsic-dichotomy}\label{fs:faithful-geometry}
Let $J=\I(a)$ have infinite rank. The following are equivalent:
\begin{enumerate}[label=\textup{(\roman*)}]
\item $J$ satisfies \eqref{fs:criterion}.
\item Whenever $\I(A)+\I(B)=J$, some $(s,t)\in\mathbb R^2\setminus\{0\}$ satisfies
$\I(sA+tB)=J$.
\item The assertion in \textup{(ii)} holds for every commuting
diagonal self-adjoint pair $A,B$.
\end{enumerate}
When these hold, every finite representation $J=\sum_i\I(A_i)$
satisfies $\I(\sum_i g_iA_i)=J$ almost surely for independent
standard real Gaussians $g_i$.
When they fail, for every $D\ge2$ there is a real-linear injection
$\Phi:\mathbb R^D\to J$, with diagonal self-adjoint values, such that
\[
 \mathcal J(U):=\sum_{z\in U}\I(\Phi(z))
       \qquad(U\subseteq\mathbb R^D\text{ a linear subspace})
\]
is principal and
\begin{equation}\label{fs:faithful-order}
 \mathcal J(U)\subseteq\mathcal J(W)\iff U\subseteq W,
 \quad \mathcal J(U+W)=\mathcal J(U)+\mathcal J(W),
 \quad \mathcal J(\mathbb R^D)=J.
\end{equation}
The map can moreover be chosen so that, simultaneously for all
$z\ne0$ and $q>0$, $\Phi(z)\in\mathcal S_q$ if and only if
$a\in\ell^q$.
\end{theorem}

For every nonnuclear $J$, the homogeneous lift $Y(x)=\Phi(x,1)$
in dimension $d+1$ therefore realizes affine geometry faithfully:
for nonempty finite $S,T\subseteq\mathbb R^d$,
\[
 \sum_{x\in S}\I(Y(x))\subseteq\sum_{t\in T}\I(Y(t))
 \iff \operatorname{aff}S\subseteq\operatorname{aff}T.
\]
Its full envelope is the prescribed $J$. Every tuple of at most
$d$ points has a strictly smaller response join, whereas an
affinely independent $(d+1)$-tuple recovers $J$.
Thus the universal $d+1$ sample bound is sharp inside every
nonnuclear envelope, including trace-class envelopes and envelopes
contained in every Schatten class.

\begin{proof}[Proof of Theorem~\ref{fs:intrinsic-dichotomy}]
Under \eqref{fs:criterion}, Lemma~\ref{fs:hausdorff} with zero
shift gives saturation for Lebesgue-almost every coefficient vector
in any finite representation of $J$. This proves the Gaussian
assertion and \textup{(i)}$\Rightarrow$\textup{(ii)}$\Rightarrow$\textup{(iii)}.
Assume now that \eqref{fs:criterion} fails. We will choose unit
vectors $u_j\in\mathbb R^D$ and set
\[
 \Phi(z)=\operatorname{diag}(a_j\langle u_j,z\rangle).
\]
For a subspace $U\subseteq\mathbb R^D$, let $P_U$ be its orthogonal
projection and put
\[
 S_U=\operatorname{diag}(a_j\|P_Uu_j\|).
\]
This is the column modulus of $\Phi$ on an orthonormal basis of $U$.
Hence $\mathcal J(U)=\I(S_U)$ is principal and
$S_{\mathbb R^D}=\operatorname{diag}(a)$, giving the prescribed top ideal.
Linearity gives monotonicity and preservation of sums.

To choose the directions, set $\delta_h=2^{-h}/(8\sqrt D)$.
At each stage $h\ge1$,
Lemma~\ref{fs:blocks}, applied to $b_j=a_{2^j}$ with power $D$
beyond all previously chosen coordinates, supplies disjoint intervals $[n_\ell,N_\ell)$ with
\[
 N_\ell>2^h n_\ell,\quad
 a_{N_\ell}\le\delta_h a_{2^h n_\ell},\quad
 \rho_\ell=\delta_h\frac{a_{2^h n_\ell}}{a_{n_\ell}},\quad
 \sum_\ell\rho_\ell^D\ge1.
\]
Rescale Lemma~\ref{fs:cover} to $[-1,1]^D$.
For each covering ball meeting the unit sphere, choose $v_\ell$
in their intersection and set $u_j=v_\ell$ on its block;
choose arbitrary unit vectors elsewhere.
Every unit vector $v$ then has a matching index with
$\|v-v_\ell\|\le4\sqrt D\rho_\ell$.

Fix $z\notin W$ and a unit vector $v\in W^\perp$ with
$\eta=|\langle v,z\rangle|>0$.
At each stage choose a block matching $v$, and abbreviate its
endpoints by $n,N$. Since $\|P_Wv_\ell\|\le4\sqrt D\rho_\ell$,
deleting the first $n-1$ coordinates gives
\[
 s_n(S_W)\le\max\{4\sqrt D\rho_\ell a_n,a_N\}
          \le4\sqrt D\delta_h a_{2^h n}.
\]
Eventually $|\langle v_\ell,z\rangle|\ge\eta/2$.
For any fixed integer $m\ge1$, once $2^h\ge m+1$, the block
contains its first $mn$ coordinates, and hence
\[
 s_{mn}(\Phi(z))\ge\frac\eta2a_{(m+1)n}
                  \ge\frac\eta2a_{2^h n}.
\]
An inclusion $\I(\Phi(z))\subseteq\I(S_W)$ would imply
$s_{mn}(\Phi(z))\le C s_n(S_W)$ for some fixed $m,C$.
The two bounds contradict this as $\delta_h\to0$.
Thus $\Phi(z)\notin\mathcal J(W)$, proving
\eqref{fs:faithful-order} simultaneously for all subspaces.
Taking $W=0$ also proves injectivity.

For $D=2$, every nonzero combination generates the ideal of a
proper line and hence a proper subideal of $J$, contradicting
\textup{(iii)}. The equivalence follows.

To retain all Schatten tests, take the divergent positive rational
exponents and the finite positive critical exponent when divergent.
Between stages add finite intervals with $u_j=e_i$ and
$\sum a_j^q\ge1$, scheduling every pair $(i,q)$ infinitely often.
Divergent tails permit this; later separating blocks begin beyond
these intervals. Any $z_i\ne0$ forces nonmembership at each tested
exponent, which by monotonicity determines all nonmembership.
Membership follows from $|\langle u_j,z\rangle|\le\|z\|$.
The unit directions preserve the top ideal and all separation estimates.
\end{proof}

\subsection{Exact source laws}
Lemma~\ref{fs:hausdorff} measures an exceptional set by the widths of
covering slabs: no cost is assigned to their tangential extent.
Let $\mathcal Z_d$ consist of the Borel sets $E\subseteq\mathbb R^d$
that, for every $s,\varepsilon>0$, admit a countable cover by open slabs
\[
 S_j=\{x:|\langle u_j,x\rangle-c_j|<r_j\},\qquad
 |u_j|=1,\quad r_j>0,\qquad \sum_jr_j^s<\varepsilon.
\]
For $d=1$ this is exactly Hausdorff dimension zero; for $d\ge2$,
every affine hyperplane belongs to $\mathcal Z_d$.
The class is preserved by countable unions and invertible affine maps.
Lemma~\ref{fs:hausdorff} puts every nuclear exceptional set in
$\mathcal Z_d$. The next theorem shows that a universal saturation law must
annihilate this class, even when the column profile is fixed.
For $d\ge2$ the covering condition is simultaneous $(s,d-1)$-tube
nullity for all $s\in(0,1]$ in the terminology of
\cite[Section~4.5]{Fraser2026}; no outer-measure regularity is used.

\begin{theorem}[The exact source condition and a universal test ideal]
\label{fs:slab-source}\label{fs:prescribed-exceptional}
\label{fs:universal-test-ideal}
Put
\[
 a_n^*=\frac1{\lfloor\log_2 n\rfloor!},\qquad
 T_*=\operatorname{diag}(a_n^*),\qquad J_*=\I(T_*).
\]
The operator $T_*$ is trace class and $J_*$ satisfies \eqref{fs:criterion}.
For a Borel probability law $\mu$ on $\mathbb R^d$ and any fixed
$\varepsilon>0$, the following are equivalent:
\begin{enumerate}[label=\textup{(\roman*)}]
\item $\mu(E)=0$ for every $E\in\mathcal Z_d$.
\item For every pencil $A(x)=\sum_{i=1}^d x_iA_i$ whose
infinite-rank coefficient ideal $J$ satisfies \eqref{fs:criterion},
and every compact $B$,
$\I(A(x)+B)=J+\I(B)$ for $\mu$-almost every $x$.
\item For every linearly independent diagonal self-adjoint list
$A_1,\ldots,A_d$ with $\sum_iA_i^2=T_*^2$ and every diagonal self-adjoint
$D\in J_*$ with $\|D\|<\varepsilon$,
$\I(A(x)+D)=J_*$ for $\mu$-almost every $x$.
\end{enumerate}
More precisely, every nonempty compact $C\in\mathcal Z_d$ is
contained in the exceptional set of a pencil and shift from
\textup{(iii)}. For $d=1$, one may always use the single fixed
operator $A(t)=tT_*$ and obtain the exact exceptional set $C$.
\end{theorem}
\begin{proof}
Lemma~\ref{fs:hausdorff} places each nuclear exceptional set in
$\mathcal Z_d$: the tails of its summable slab covers have arbitrarily
small cost. This proves \textup{(i)}$\Rightarrow$\textup{(ii)}.
For $2^j\le n<2^{j+1}$,
$a_{4n}^*/a_n^*=1/((j+1)(j+2))$, so
\[
 \sum_{n\ge1}\frac{a_{4n}^*}{na_n^*}
 \le\sum_{j\ge0}\frac1{(j+1)(j+2)}<\infty,
 \qquad \operatorname{tr}T_*=\sum_{j\ge0}\frac{2^j}{j!}=e^2.
\]
Thus \textup{(ii)} implies \textup{(iii)}.

To prove \textup{(iii)}$\Rightarrow$\textup{(i)}, it suffices by
inner regularity to place every nonempty compact $C\in\mathcal Z_d$
inside an exceptional set from \textup{(iii)}.
The factorial profile converts a small $\sum_i r_i^{1/h}$ cover
into the required block ratios. Fix such a $C$ and put
$R=\sup_{x\in C}\|x\|$, $M=\max\{1,2R\}$.
Choose $n_0>d$ with $R a_{n_0}^*<\varepsilon$; all blocks below
will start after this initial segment.
At stage $h\ge1$, put $\delta=2^{-h}$ and choose $j_0\ge1$ with
$2^{j_0}\ge n_0$, beyond all earlier blocks, and
$j_0+h+1\ge M/\delta$.
Take a finite open slab cover
$|\langle u_i,x\rangle-c_i|<r_i$, with $\|u_i\|=1$,
$|c_i|\le R$, and decreasingly ordered halfwidths satisfying
\[
 \sum_i r_i^{1/h}\le\frac{\delta^{1/h}}{j_0+h}.
\]
Compactness and $C\in\mathcal Z_d$ supply such a cover.
Its centers can be chosen by moving each middle hyperplane through
a covered point of $C$ and doubling its halfwidth, with the initial
budget reduced accordingly.
Set $j_i=j_0+(h+1)(i-1)$, $n_i=2^{j_i}$ and
$N_i=2^{j_i+h+1}$. Since $j_i+h\le(j_0+h)i$ and
$i r_i^{1/h}\le\sum_\ell r_\ell^{1/h}$,
\[
 r_i\le\frac\delta{(j_i+h)^h}
       \le\delta\frac{j_i!}{(j_i+h)!}
       =\delta\frac{a_{2^h n_i}^*}{a_{n_i}^*},\qquad
 \frac{a_{N_i}^*}{a_{2^h n_i}^*}
       =\frac1{j_i+h+1}\le\frac\delta M.
\]
The blocks $[n_i,N_i)$ are disjoint. Assign their normals and
centers to their coordinates; on gaps choose unit normals and
centers in $[-R,R]$, taking $c_n=0$ for $n<n_0$.
Write these data as $u_n,c_n$ and define
\[
 A_i=\operatorname{diag}(a_n^*u_{n,i}),\qquad
 D=\operatorname{diag}(-a_n^*c_n).
\]
Then $\|D\|\le R a_{n_0}^*<\varepsilon$.
Since $\|u_n\|=1$, the column modulus is $T_*$;
hence $\sum_i\I(A_i)=J_*$ and $D\in J_*$.
Choose the first $d$ normals to be the coordinate basis,
ensuring linear independence without changing any selected block.

For $x\in C$, select a covering slab at each stage. Deleting
the first $n_i-1$ coordinates yields
\[
 s_{n_i}(A(x)+D)
 \le\max\{r_i a_{n_i}^*,Ma_{N_i}^*\}
 \le\delta a_{2^h n_i}^*.
\]
As $h\to\infty$, this contradicts every finite-dilation comparison
required for $J_*\subseteq\I(A(x)+D)$.
For $d=1$, orient every interval with normal $+1$ and choose a
point of $C$ as the center on each gap after the initial segment.
Then $A(t)=tT_*$ and all tail centers lie in $C$.
For $t\notin C$, the tail multiplier $t-c_n$ is boundedly invertible.
Finite initial segments do not change its infinite-rank ideal,
so $\I(tT_*+D)=J_*$, proving exact scalar realization.
\end{proof}

\begin{remark}[No countable collection of test pencils]
\label{fs:no-countable-test}
For $d\ge2$, the pencil quantifier in Theorem~\ref{fs:slab-source}
cannot be replaced by a countable list, even if every compact shift
is tested. Indeed, let $A^{(m)}$ be such a list, with infinite-rank
nuclear coefficient ideals $J_m$.
Theorem~\ref{fs:main} at shift zero gives a common nonzero direction
$v$ with $\I(A^{(m)}(v))=J_m$ for every $m$.
Let $\mu$ be uniform probability on $C=\{tv:0\le t\le1\}$.
For each $m$ and fixed compact $B$, the scalar case of that theorem gives
\[
 \I(A^{(m)}(x)+B)=J_m+\I(B)
       \qquad\text{for }\mu\text{-almost every }x.
\]
Nevertheless $C\in\mathcal Z_d$ and $\mu(C)=1$.
Thus one line-supported source passes all these tests while violating
the source condition. In particular, even with the column modulus
fixed as $T_*$, the realizing coefficients must depend on the prescribed set.
\end{remark}

\begin{remark}[Affine sections constrain exact prescription]
\label{fs:affine-flat-boundary}
Suppose $B\in J$ and $J$ satisfies \eqref{fs:criterion}.
For every $k$-dimensional affine flat $L$, $k\ge1$, either
$L\subseteq E_B$ or $\dim_{\mathrm H}(E_B\cap L)\le k-1$.
Indeed, let $F_L$ be the ideal generated by $A(x)+B$, $x\in L$.
If $F_L\subsetneq J$, the whole flat is exceptional.
Otherwise, homogenizing the restricted affine pencil gives
$k+1$ coefficients with ideal $J$; Lemma~\ref{fs:hausdorff}
bounds its exceptional cone by dimension $k$.
On bounded charts, $(t,y)\mapsto(t,ty)$, $1\le t\le2$, is
bi-Lipschitz and sends $[1,2]\times(E_B\cap L)$ into that cone.
Hausdorff slicing gives the asserted bound.
Thus a nondegenerate proper line segment cannot be the exact
$E_B$, although for $d\ge2$ it lies in an affine hyperplane.
\end{remark}

For any fixed probability law $\mu$ annihilating $\mathcal Z_d$,
condition \textup{(v)} of Theorem~\ref{fs:main} remains equivalent
when the Gaussian vectors are replaced by independent
$\xi,\xi'\sim\mu$. Indeed, nuclearity gives $\mu(E_B)=0$ for every
fixed $B$. Conversely, hyperplanes belong to $\mathcal Z_d$,
so in the nonnuclear case the same residual set of shifts gives
almost-sure incomparability of the two responses.

A sufficient condition for a probability law $\mu$ to annihilate
$\mathcal Z_d$ is the uniform slab bound
\begin{equation}\label{fs:nonconcentration}
 \sup_{|u|=1,\ c\in\mathbb R}
 \mu\{x:|u\cdot x-c|\le t\}\le C_\mu t^\alpha
                         \qquad(0<t\le1)
\end{equation}
for some $C_\mu<\infty$ and $\alpha>0$.
Summing over covers proves the assertion; increasing $C_\mu$
extends the bound to all $t>0$.

For example, uniform probability on $S^{d-1}$, $d\ge2$, satisfies
such a bound with exponent $1/2$ for $d=2$ and $1$ for $d\ge3$:
its one-coordinate density is proportional to
$(1-t^2)^{(d-3)/2}$. Thus annihilating every set of ordinary
Hausdorff dimension at most $d-1$ would be unnecessarily restrictive.

\begin{proposition}[Absolute continuity of all projections is insufficient]
\label{fs:projection-boundary}
For every $d\ge2$, there is a compactly supported probability measure
$\mu$ on $\mathbb R^d$ whose every nonzero one-dimensional linear
projection is absolutely continuous, although $\mu(S)=1$ for a
compact set $S\in\mathcal Z_d$.
\end{proposition}
\begin{proof}
Choose a perfect nowhere-dense set $F\subset[0,1]$ of positive
Lebesgue measure, and a smooth function $g\ge0$ vanishing exactly
on $F$. Such a function is obtained by summing positive smooth bumps
on the complementary intervals with sufficiently small coefficients.
Set
\[
 v(t)=\int_0^t g(u)\,du,\qquad
 \gamma(t)=\int_0^t(1,v(u),\ldots,v(u)^{d-1})\,du,
\]
and let $\mu$ be the image of normalized Lebesgue measure on $F$
under $\gamma$. Since $F$ has empty interior, $v$ is strictly
increasing. For every nonzero $\theta\in\mathbb R^d$,
$(\theta\cdot\gamma)'$ is a nonzero polynomial in $v$ and has
at most $d-1$ zeros. Away from these points the projection is locally
a $C^1$ diffeomorphism, so its pushforward measure is absolutely
continuous.

Perfectness implies that every derivative of $g$ vanishes on $F$,
and hence $\gamma^{(k)}=0$ there for every $k\ge2$.
For each integer $N\ge2$, Taylor's theorem shows that the part of
$\gamma(F)$ within parameter distance $r$ of $t\in F$ lies in a
slab through $\gamma(t)$, containing its tangent line, of width
$O_N(r^N)$. Covering $F$ by $O(r^{-1})$ such intervals gives,
for every $s>0$, a slab-cover cost $O_N(r^{Ns-1})$.
Choose $N>1/s$ and let $r\downarrow0$. Thus
$S=\gamma(F)$ belongs to $\mathcal Z_d$.
\end{proof}

In dimension two, this example also shows that convolution can restore
complete determinacy. For fixed $a,b\ne0$, the map
$(t,s)\mapsto a\gamma(t)+b\gamma(s)$ has Jacobian
\[
 ab\bigl(v(s)-v(t)\bigr)\ne0\qquad(t\ne s).
\]
Thus $a\xi+b\xi'$ has an absolutely continuous law for independent
$\xi,\xi'\sim\mu$, by countably many local inverse-function charts.
Use Theorem~\ref{fs:slab-source} to choose a nuclear diagonal pencil
$A$, with coefficient ideal $J$, and $D\in J$ whose exceptional
set contains $S$. Then $X=A(\xi)+D$ is bounded in trace norm,
diagonal and self-adjoint. Its unshifted ideal has no deterministic
principal-ideal atom by Proposition~\ref{fs:response-defects}, whereas
Theorem~\ref{fs:main} gives, for every fixed $a,b\ne0$ and compact $B$,
\[
 \I(aX+bX'+B)=J+\I(B)\qquad\text{almost surely}.
\]
Thus complete determinacy of the convolution square or symmetrization
does not imply complete determinacy of $X$. The coupled absorption
condition in Theorem~\ref{ideal:response-structure} is substantive.

\subsection{Absolutely continuous polynomial sources}

Let $Z$ be an $\mathbb R^k$-valued random vector with an arbitrary
absolutely continuous law. Expand a finite compact-operator-valued
polynomial of degree at most $r\ge1$ as
\begin{equation}\label{fs:model}
 X=P(Z)=A_0+\sum_{i=1}^d H_i(Z)A_i,
 \qquad J=\I(a),\qquad
 a_n=s_n\Bigl(\bigl(\sum_iA_i^*A_i\bigr)^{1/2}\Bigr).
\end{equation}
Here the $H_i$ are distinct nonconstant real Hermite polynomials,
normalized in standard Gaussian $L^2$, and zero operator coefficients
are omitted. The operators $A_i$ need not be linearly independent.
This is an algebraic choice of polynomial basis; no moment of $Z$
is assumed. The triangular change from monomials preserves the
nonconstant coefficient ideal, and $A_0-P(0)\in J$. Thus $J$ and
the coset $-A_0+J$ are intrinsic to $P$.
Every nontrivial affine combination of the $H_i$ is a nonzero
polynomial, so their coefficient law gives affine hyperplanes zero
mass. Proposition~\ref{fs:response-defects} therefore applies even
when that law has no density in $\mathbb R^d$.

\begin{corollary}[Polynomial sources with a density]\label{fs:polynomial}
Assume $J$ in \eqref{fs:model} contains an infinite-rank operator.
For any prescribed absolutely continuous law of $Z$, the following
are equivalent: \eqref{fs:criterion} holds; and for every fixed
compact $B$, $\I(P(Z)+B)$ is almost surely deterministic.
Under these equivalent conditions,
\begin{equation}\label{fs:response}
 \I(P(Z)+B)=J+\I(A_0+B)=J+\I(P(0)+B)
 \qquad\text{almost surely}
\end{equation}
for each fixed compact $B$. If the condition fails, there is an
arbitrarily small $D\in J$ such that $\I(P(Z)-A_0+D)$ is
nondeterministic. For a residual set of compact shifts $B$, two
independent responses at $B$ are almost surely incomparable.
The finite-sample and order conclusions of
Proposition~\ref{fs:response-defects} apply in either case.
\end{corollary}
\begin{proof}
First evaluate the same polynomial at a standard Gaussian $G$ in
$\mathbb R^k$ and put $\xi=(H_1(G),\ldots,H_d(G))$.
For $|u|=1$ and $c\in\mathbb R$, the polynomial
$Q=u\cdot\xi-c$ has degree at most $r$ and
$\|Q\|_2^2=1+c^2\ge1$. Carbery--Wright
\cite[Theorem~8, with parameter $2r$]{CW} gives
\begin{equation}\label{fs:polynomial-slab}
 \Pbb\{|u\cdot\xi-c|\le t\}
 \le C_0r\,t^{1/r}\|Q\|_2^{-1/r}
 \le C_0r\,t^{1/r}\qquad(t>0).
\end{equation}
Summing over slab covers shows that the law of $\xi$ annihilates
$\mathcal Z_d$. Under \eqref{fs:criterion},
Theorem~\ref{fs:slab-source}, with shift $A_0+B$, gives
\eqref{fs:response} for $P(G)$.
For each fixed $B$, its exceptional set in $\mathbb R^k$ is Borel
by the finite-dilation comparisons. Gaussian measure and Lebesgue
measure have the same null sets, so absolute continuity transfers
the identity to $Z$. The expression using $P(0)$ follows from
$A_0-P(0)\in J$.

Conversely, the law of $(H_1(Z),\ldots,H_d(Z))$ gives every affine
hyperplane zero mass. If \eqref{fs:criterion} fails,
Lemma~\ref{fs:failure} supplies the small internal $D$ and proves
nondeterminacy. Proposition~\ref{fs:generic-failure} gives the
residual assertion, since translation by $A_0$ preserves residual
sets. The remaining conclusions follow from
Proposition~\ref{fs:response-defects}.
\end{proof}

This proof transfers only a null set from a reference Gaussian law.
The actual source $Z$ needs no moments or independent coordinates;
independent Cauchy coordinates, for example, are allowed.
The affine description of $E_B$ following Lemma~\ref{fs:hausdorff},
together with hyperplane-nullity of the coefficient law, also gives
\eqref{fs:response} for nonzero finite-rank $J$ by using the shift
$A_0+B$. If $J=0$, the model is deterministic.

We now use the rate in \eqref{fs:nonconcentration} to control
spectral comparison constants. The Gaussian coefficient law in
\eqref{fs:polynomial-slab} satisfies this bound with $\alpha=1/r$.

\begin{lemma}[Uniform lower spectral control]\label{fs:lower}
For the pencil $A(x)=\sum_i x_iA_i$, suppose that $J$ has infinite
rank and that \eqref{fs:criterion} and \eqref{fs:nonconcentration}
hold. Let $\xi\sim\mu$. There are a deterministic integer $M_0$
and a constant $C_1$, independent of the compact shift, such that
for each fixed compact $B$ there is a random $0<c_B\le1$ with
\begin{equation}\label{fs:lower-tail}
 s_n(A(\xi)+B)\ge c_Ba_{M_0n}\quad(n\ge1),\qquad
 \sup_{B\in\mathcal K(H)}\Pbb\{c_B<\varepsilon\}
 \le C_1\varepsilon^\alpha\quad(0<\varepsilon<1).
\end{equation}
The spectral inequalities hold simultaneously for all $n$ on an
event of probability one that may depend on $B$. Moreover,
$\sup_B\E c_B^{-p}<\infty$ for $0<p<\alpha$, and
$\I(A(\xi)+B)=J+\I(B)$ almost surely for each fixed $B$.
\end{lemma}
\begin{proof}
Choose an integer $C$ with $2^C\ge2d$. By
Lemma~\ref{fs:normalization}, there is $K$ for which
\[
 C_\alpha=\sum_{j\ge0}
 \left(\frac{a_{2^{j+C+K}}}{a_{2^{j+C}}}\right)^\alpha<\infty.
\]
For each fixed $B$, define the measurable variable
\[
 c_B=\min\left\{1,\inf_{j\ge0}
       \frac{s_{2^j}(A(\xi)+B)}{a_{2^{j+C+K}}}\right\}.
\]
Lemma~\ref{fs:slab}, \eqref{fs:nonconcentration}, and a union bound give
\[
 \Pbb\{c_B<\varepsilon\}
 \le C_2\sum_{j\ge0}
 \left(\frac{\varepsilon a_{2^{j+C+K}}}{a_{2d2^j}}\right)^\alpha
 \le C_2C_\alpha\varepsilon^\alpha.
\]
Letting $\varepsilon\downarrow0$ gives $c_B>0$ almost surely.
For $2^j\le n<2^{j+1}$, monotonicity yields
\[
 s_n(A(\xi)+B)\ge c_Ba_{2^{j+1+C+K}}
                 \ge c_Ba_{2^{C+K+1}n}.
\]
This proves \eqref{fs:lower-tail} with $M_0=2^{C+K+1}$.
Tail integration gives the inverse moments.
The lower bound puts $J$ in the response ideal; subtraction then
puts $B$ there, proving the asserted identity.
\end{proof}

\begin{corollary}[Gaussian full-profile comparison]\label{fs:gaussian}
Suppose the infinite-rank ideal $J$ in \eqref{fs:model} satisfies
\eqref{fs:criterion} and $Z=G$ is standard Gaussian. Set
$b_n(B)=\max\{a_n,s_n(A_0+B)\}$.
There is a deterministic integer $M\ge1$, independent of $B$,
such that for each fixed compact $B$ a random $L_B\ge1$ satisfies,
on one event of probability one for all $n\ge1$,
\begin{equation}\label{fs:profile}
 L_B^{-1}b_{Mn}(B)\le s_n(P(G)+B)
       \le L_Bb_{\lceil n/M\rceil}(B),\qquad
 \sup_{B\in\mathcal K(H)}\E L_B^p<\infty\quad(0<p<1/r).
\end{equation}
The moment range is optimal over degree-$r$ models.
\end{corollary}
\begin{proof}
By \eqref{fs:polynomial-slab}, Lemma~\ref{fs:lower} applies to
$\xi=(H_1(G),\ldots,H_d(G))$ with $\alpha=1/r$.
Write $F=A_0+B$, $T=A(\xi)+F$, $\ell=|\xi|$, and
$b_n=\max\{a_n,s_n(F)\}$.
The lemma gives $a_{M_0n}\le c_F^{-1}s_n(T)$, while column
factorization gives $s_n(A(\xi))\le\ell a_n$.
The singular-value sum inequality therefore yields
\[
 s_n(T)\le(1+\ell)b_{\lceil n/2\rceil},\qquad
 s_{(M_0+1)n}(F)\le(1+\ell/c_F)s_n(T).
\]
Combining the second bound with that for $a_{M_0n}$ gives
\[
 b_{(M_0+1)n}\le\left(1+\frac{\ell+1}{c_F}\right)s_n(T).
\]
Take $M=M_0+1$ and $L_B=1+(\ell+1)/c_F$.
Since $M\ge2$ and $L_B\ge1+\ell$, monotonicity of $b$ also gives
the upper comparison in \eqref{fs:profile}.
All positive moments of $\ell$ are finite. Given $p<1/r$, choose
$q>1$ with $pq<1/r$ and apply H\"older to
$(\ell+1)^pc_F^{-p}$. Lemma~\ref{fs:lower} bounds its expectation
uniformly in $F$, proving the moment assertion without any
independence assumption between $\ell$ and $c_F$.

For optimality, take $X=g^rT$ with $g$ standard Gaussian and
$s_n(T)=e^{-n}$. The coefficient ideal satisfies the criterion.
For every fixed finite $M$, the lower comparison at $n=1$ forces
$L_0\ge c_{r,M}|g|^{-r}$ with $c_{r,M}>0$.
Hence $\E L_0^p=\infty$ for $p\ge1/r$.
\end{proof}

Absolute continuity alone gives no positive moments for the
comparison constant. If $\xi$ has density
$[t(\log(1/t))^2]^{-1}$ on $(0,e^{-1})$ and
$T=\operatorname{diag}(e^{-n})$, the qualitative result applies to
$X=\xi T$. But every fixed-dilation lower comparison at $n=1$
forces $L_0\ge c_M/\xi$, so $\E L_0^p=\infty$ for every $p>0$.

\section{Conditional contraction and spectral saturation}
\label{cov:section}\label{ops:section}\label{linear:section}
Call a centered countable array $Z$ a \emph{curvature array} if every
finite coordinate marginal has positive density $ce^{-V}$ and
\begin{equation}\tag{H$_0$}\label{ext:array-curvature}
 \alpha I\le\nabla^2V\le\beta I,\qquad 0<\alpha\le\beta<\infty,
\end{equation}
with constants independent of the finite set.
These bounds concern joint coordinate laws and allow dependence
between rows. Potentials may be $C^{1,1}_{\mathrm{loc}}$, with the
inequalities holding almost everywhere.
To track linear images, we use a Gaussian reference with covariance
$C>0$ on its support $F$, and write the corresponding condition as
\begin{equation}\tag{H}\label{curv:class}
 \alpha C^{-1}\le\nabla_F^2V\le\beta C^{-1}.
\end{equation}

All covariance inverses are taken on their supports.
The lower Hessian bound gives the moment and convex-order estimates
used for upper spectral bounds. The upper Hessian bound gives a
shift-uniform Laplace estimate, hence lower spectral bounds.
Stability under linear images and centered conditional limits lets us
apply both estimates after conditioning. The square-tail product
records the generator obtained by this iteration.

\subsection{Curvature on conditional fibers}

\begin{lemma}[Conditional covariance and Laplace bounds]\label{var:laplace}
Let $W$ have positive density $ce^{-V}$ on a finite-dimensional support $F$,
and let $C>0$ on $F$ satisfy \eqref{curv:class}.
For linear observation and target maps $O,A$, and $g\sim N_F(0,C)$, put
$S=\operatorname{Cov}(Ag-\E[Ag\mid Og])$. Almost surely,
\begin{equation}\label{curv:fiber-covariance}
 \beta^{-1}S\le\operatorname{Cov}(AW\mid OW)\le\alpha^{-1}S.
\end{equation}
For an observation-measurable pair $(z,P)$, with $P$ an orthogonal projection,
and $t>0$,
\begin{equation}\label{var:laplace-order}
 \E[e^{-t\|P(AW-z)\|^2/2}\mid OW]
 \le\det(I+(t/\beta)PSP)^{-1/2}.
\end{equation}
If $W$ is centered, every nonnegative convex test satisfies
\begin{equation}\label{curv:convex-order}
 \E\Phi(W)\le\E\Phi(\alpha^{-1/2}g).
\end{equation}

\end{lemma}
\begin{proof}
Work inside $F$ and choose an isometry $J$ onto $\ker O$.
The conditional fiber is $w_y+Jv$, with potential $V(w_y+Jv)$; the coarea
Jacobian is constant on that fiber. If $H=J^*C^{-1}J$, its Hessian lies between
$\alpha H$ and $\beta H$, while completion of the Gaussian square gives
\begin{equation}\label{curv:fiber-identity}
 S=AJH^{-1}J^*A^*.
\end{equation}
A zero-dimensional fiber is deterministic. No inverse of a singular ambient
covariance has been used.

For a smooth potential $U$ with $aI\le\nabla^2U\le D_0$, integration by parts gives
\[
 \E\nabla U=0,\quad
 \E[(Z-\E Z)(\nabla U)^*]=I,\quad
 \E[(\nabla U)(\nabla U)^*]=\E\nabla^2U.
\]
The block covariance Schur complement yields the Cram\'er--Rao bound
$\operatorname{Cov}(Z)\ge(\E\nabla^2U)^{-1}\ge D_0^{-1}$.
For $\varphi(t)=\E e^{-t\|DZ-z\|^2/2}$, let $\nu_t$ have potential
$U_t(x)=U(x)+t\|Dx-z\|^2/2$, and put
$Q_t=(D_0+tD^*D)^{-1}$.
The same covariance bound gives $\operatorname{Cov}_{\nu_t}(Z)\ge Q_t$.
Jacobi's determinant formula yields
\[
\begin{split}
 -\frac d{dt}\log\!\left[
   \varphi(t)\det(I+tDD_0^{-1}D^*)^{1/2}\right]
 &=\tfrac12\operatorname{tr}
       \bigl[D(\operatorname{Cov}_{\nu_t}(Z)-Q_t)D^*\bigr]\\
 &\quad+\tfrac12\|D\E_{\nu_t}Z-z\|^2\ge0.
\end{split}
\]
Since $\varphi(0)=1$, this proves
$\varphi(t)\le\det(I+tDD_0^{-1}D^*)^{-1/2}$ uniformly in $z$.
Strong convexity supplies Gaussian tails and the upper Hessian bound gives
linear gradient growth, justifying the integrations by cutoff.
Use $D_0=\beta H$ and $D=PAJ$ on the fiber to obtain
\eqref{var:laplace-order} and the lower half of \eqref{curv:fiber-covariance}.
Brascamp--Lieb's variance inequality \cite[Theorem~4.1]{BL1976}, applied to fiberwise
linear functions, gives the upper half.

The density of centered $W$ relative to $N_F(0,\alpha^{-1}C)$ is log-concave;
Harg\'e's theorem \cite[Theorem~1.1]{Harge2004} gives \eqref{curv:convex-order}.
These statements extend to $C^{1,1}_{\mathrm{loc}}$ potentials by convolution
smoothing of the potential: both Hessian bounds persist, the normalizations
converge with Gaussian domination, and centering the approximants gives the
convex comparison. Fiberwise bounds hold for all centers before passage to the
limit.
\end{proof}

\begin{lemma}[Linear images and weak limits]\label{curv:stability}
A centered finite-dimensional law satisfying \eqref{curv:class} relative to
$C$ retains that condition under a linear map $L$, relative to $LCL^*$
on its range, with the same $\alpha,\beta$.
On a fixed finite-dimensional space $F$, suppose centered laws $\mu_j$
with positive densities satisfy the same bounds relative to $C_j>0$.
If $C_j\to C>0$ in operator norm and $\mu_j\Rightarrow\mu$, then $\mu$
is centered and has a positive density with a $C^{1,1}_{\mathrm{loc}}$
potential satisfying those bounds relative to $C$.
This also applies to tests with $0\le C_j\uparrow C$: on
$F=\operatorname{ran}C$, the reference covariances have full rank eventually.
\end{lemma}
\begin{proof}
For the linear-image assertion, whiten on the original reference support
and split into orthogonal coordinates $(x,y)$ so that the map depends only
on $x$. The marginal potential obeys
\[
 W_{xx}=\E[V_{xx}\mid x]-\operatorname{Cov}(V_x\mid x),\qquad
 \operatorname{Cov}(V_x\mid x)
 \le\E[V_{xy}V_{yy}^{-1}V_{yx}\mid x].
\]
The second inequality is Brascamp--Lieb on the $y$-fiber. The Schur
complement of $\nabla^2V\ge\alpha I$ gives $W_{xx}\ge\alpha I$;
the first identity gives $W_{xx}\le\beta I$. Transform back on
$\operatorname{ran}(LCL^*)$. Smooth approximation treats $C^{1,1}$ potentials.

For weak closure, discard finitely many indices so $C_j$ and $C_j^{-1}$
are uniformly bounded. Normalize $V_j(0)=0$. Integration by parts and
convex comparison give
\[
 \|\nabla V_j(0)\|
 \le\beta\|C_j^{-1}\|\E_{\mu_j}\|x\|\le M<\infty.
\]
The Hessian bounds now give locally equibounded, equicontinuous potentials
and gradients, and a common lower bound $V_j(x)\ge c\|x\|^2-M\|x\|$.
Arzel\`a--Ascoli gives subsequential local $C^1$ convergence. Quadratic
domination passes the normalizations and first moments to the limit,
identifying a positive centered density. Second-difference inequalities
pass both Hessian bounds to its $C^{1,1}_{\mathrm{loc}}$ potential.
Uniqueness of the weak limit identifies the law. If $C_j\uparrow C$,
restrict to $\operatorname{ran}C$; finite-dimensional eigenvalue convergence
makes $C_j$ positive there eventually. No singular ambient inverse is used.
\end{proof}

\subsection{The square-tail product}
Write $F_u(n)=\sum_{k\ge n}u_k^2$.
For nonnegative nonincreasing $u,v\in\ell^2$, define their
square-tail product by grouping the weights $u_i^2v_j^2$
according to the smaller index:
\begin{equation}\label{ext:square-tail-product}
 (u\odot v)_n^2
 =\sum_{\min\{i,j\}=n}u_i^2v_j^2
 =u_n^2F_v(n)+v_n^2F_u(n+1).
\end{equation}
The positive expansion shows that the sequence is nonincreasing;
its support length is the smaller of those of $u,v$.
Summing over indices at least $n$ gives
\begin{equation}\label{ext:tail-identity}
 F_{u\odot v}(n)=F_u(n)F_v(n),\qquad
 \|u\odot v\|_2=\|u\|_2\|v\|_2.
\end{equation}
Square tails determine their sequence, so $\odot$ is commutative
and associative.
For $(\Delta_d u)_n=u_{\lceil n/d\rceil}$ and
$k=\lceil n/d\rceil$, $F_{\Delta_d u}(n)\le dF_u(k)$; hence
\begin{equation}\label{ext:dilation-product}
 \Delta_d u\odot\Delta_d v\le\sqrt{2d}\,\Delta_d(u\odot v).
\end{equation}
Indeed its squared left side is at most
$d(u_k^2F_v(k)+v_k^2F_u(k))\le2d(u\odot v)_k^2$.
If $u\le C\Delta_d w$, monotonicity and $v\le\Delta_d v$ imply
$u\odot v\le C\sqrt{2d}\,\Delta_d(w\odot v)$.
Thus either input may be replaced by any generator of the same principal
ideal. This compatibility is what permits conditional iteration.

\subsection{Conditional distances}
For lower singular-value bounds after a shift, we condition on the
other rows and estimate a row's distance from their shifted span.
The Laplace bound controls its inverse moments; the inverse Gram
identity then transfers them to singular values.
\begin{lemma}[Inverse moments for tensor rows]\label{var:lowerprob}
Let $Z$ be a curvature array, let $a_1,\ldots,a_m>0$, and let
$b\in\ell^2$ be nonnegative and nonincreasing.
Let $Y:\ell^2\to\mathbb R^m$ have entries $Y_{ik}=a_iZ_{ik}b_k$.
For every fixed bounded $F:\ell^2\to\mathbb R^m$, every integer $p\ge1$
with $F_b(m+2p)>0$, and $1\le j\le m$,
\begin{equation}\label{curv:negative-moment}
 \|s_j(Y+F)^{-2}\|_{L^p}
 \le \frac{\beta}{m-j+1}
       \sum_{i=1}^m\frac{1}{a_i^2F_b(m+2p)}.
\end{equation}
\end{lemma}
\begin{proof}
First retain only finitely many input coordinates. Given all other
rows, project off their shifted span, and let $d_i$ be the distance
of the $i$th shifted row from that span. The projection is measurable
with respect to the other rows and has codimension at most $m-1$.
The Gaussian reference rows are independent; their residual
covariances are $a_i^2\operatorname{diag}(b_k^2)$.
Lemma~\ref{var:laplace} and min--max therefore give
\[
 \E[e^{-td_i^2/2}\mid\text{other rows}]
 \le \prod_{k\ge m}
       (1+(t/\beta)a_i^2b_k^2)^{-1/2}.
\]
Take the cutoff large enough that its zero-completed tail at
$m+2p$ is positive. For this finite list put
$q=\beta^{-1}a_i^2b_{m+2p}^2$ and
$\tau=\beta^{-1}a_i^2F_b(m+2p)$.
The first $2p$ factors have eigenvalues at least $q$.
For the remaining eigenvalues $0\le x\le q$,
$\log(1+tx)\ge(x/q)\log(1+tq)$.
The displayed bound is thus at most $(1+tq)^{-p-\tau/(2q)}$.
It tends to zero as $t\to\infty$, so $d_i>0$ almost surely.
The integral representation of a negative power gives
\[
 \E[d_i^{-2p}\mid\text{other rows}]
 \le\frac{1}{2^p\Gamma(p)}\int_0^\infty
       t^{p-1}(1+tq)^{-p-\tau/(2q)}\,dt
 =\prod_{\ell=0}^{p-1}(\tau+2\ell q)^{-1}\le\tau^{-p}.
\]
The inverse Gram identity gives
\[
 \sum_{\ell=1}^m s_\ell(Y+F)^{-2}=\sum_{i=1}^m d_i^{-2}.
\]
Its last $m-j+1$ terms are at least
$(m-j+1)s_j(Y+F)^{-2}$; Minkowski proves
\eqref{curv:negative-moment} for the finite input.
For the full input, apply this argument to the first $N$ columns.
The finite-dimensional codomain makes $F$ Hilbert--Schmidt.
The corresponding tail sums tend to $F_b(m+2p)$, while
$(Y+F)Q_N\to Y+F$ almost surely in Hilbert--Schmidt norm.
Here $Y$ is Hilbert--Schmidt almost surely because
$\E\|Y\|_2^2\le\alpha^{-1}\sum_i a_i^2\|b\|_2^2$.
Singular-value continuity and Fatou give the result.
\end{proof}

\subsection{Linear contraction and conditional innovation}
\begin{lemma}[Linear contraction and conditional innovation]\label{ext:curved-linear}
For a curvature array $Z$, deterministic Hilbert--Schmidt $L,R$ and fixed
compact $F$, the linear series $LZR$ belongs to $L^2(\Omega;\mathcal S_2)$ and
\begin{equation}\label{ext:curved-linear-law}
 \E\|LZR\|_2^2\le\alpha^{-1}\|L\|_2^2\|R\|_2^2,\qquad
 \I(LZR+F)=\I(s(L)\odot s(R))+\I(F)\quad\text{a.s.}
\end{equation}
Now let $(U,V)$ be disjoint coordinate subarrays of one curvature array,
$u,v\in\ell^2$ nonnegative and nonincreasing, and $K\in\mathcal S_2$.
Suppose the coordinate truncations of $D_uUKV^*D_v$ along a fixed
exhaustion, possibly after subtraction of deterministic Hilbert--Schmidt
operators $C_N$, converge to $B$ in $L^2(\Omega;\mathcal S_2)$. The subtractions need not converge;
$C_N=0$ is the raw case. With $\mathcal F=\sigma(V)$ and
$M=\E[B\mid\mathcal F]$, every fixed compact $F$ satisfies
\begin{equation}\label{ext:prediction-law}
 \I(B+F)=\I(u\odot s(K)\odot v)+\I(M+F)\quad\text{a.s.}
\end{equation}
No conditional-centering assumption on $U$ is made.
\end{lemma}
\begin{proof}
\emph{Construction and active coordinates.}
Covariance domination gives
$\E\|LZR\|_2^2\le\alpha^{-1}\|L\|_2^2\|R\|_2^2$
for finite coefficients. Applied to coefficient differences, it
constructs the claimed $L^2(\mathcal S_2)$ limit.
For the ideal identity, it suffices by Lemma~\ref{ideal:completion}
to prove upper membership and lower inclusion for every fixed shift.
Singular decompositions of $L,R$ reduce these two tasks to the active
input and output supports, with the shift compressed to those supports.
Bounded multiplication preserves upper membership, and compression
can only decrease the ideal of the shifted operator; completion then
recovers the response to the original shift. On the active supports
write $T=D_aZD_b$, where $a=s(L)$ and $b=s(R)$; the transformed
coordinate array retains \eqref{ext:array-curvature} by
Lemma~\ref{curv:stability}.
Indeed each transformed coordinate is an $\ell^2$ linear functional
of the original array, and the coefficient vectors of any finite list
are orthonormal. Finite coordinate truncations have Gaussian reference
covariances increasing to the identity. Covariance domination gives
their $L^2$ convergence, and the weak-closure lemma preserves
\eqref{ext:array-curvature}. This uses no samplewise boundedness of $Z$.
Put $h=a\odot b$ and
\[
 M_n=\max\{a_n^2F_b(n),\,b_n^2F_a(n)\}.
\]
The square-tail formula gives $M_n\le h_n^2\le2M_n$.

\emph{Upper membership.}
For independent standard Gaussian entries $G$, Gaussian comparison
and concentration \cite[Theorems 6.8 and 3.25]{vH} give, for $p\ge1$,
\[
 \|D_cGD_d\|_{L^p({\rm op})}
 \le C\bigl(\|c\|_\infty\|d\|_2+\|c\|_2\|d\|_\infty
                 +\sqrt p\,\|c\|_\infty\|d\|_\infty\bigr).
\]
Indeed comparison of bilinear increments bounds the mean by the
first two terms, and the weak variance is at most
$\|c\|_\infty^2\|d\|_\infty^2$.
Convex comparison \eqref{curv:convex-order} transfers this bound
to $Z$ with the factor $\alpha^{-1/2}$.
Delete the first $2n-2$ rows and columns and take $p=n$.
The deleted rows and columns each cost rank at most $2n-2$.
Moreover $n a_{2n-1}^2\le F_a(n)$ and
$n b_{2n-1}^2\le F_b(n)$.
Consequently
\begin{equation}\label{lin:upper-budget}
 \|s_{4n-3}(T)\|_{L^n}\le C\alpha^{-1/2}h_n.
\end{equation}
Finite compressions and Fatou justify the bound for infinite
coefficients. If $h_n=0$, this moment bound directly gives
$s_{4n-3}(T)=0$ almost surely. At the positive terms, Markov and
Borel--Cantelli imply $T\in\I(h)$ almost surely.

\emph{Lower inclusion under fixed shifts.}
Take $n\ge8$, set $j=\lfloor n/2\rfloor$,
$m=\lfloor(n+j)/2\rfloor$, and
$p=\lfloor(n-j+1)/4\rfloor$. Then $m+2p\le n$.
Compress $T+F$ to the first $m$ rows.
Lemma~\ref{var:lowerprob}, $a_i\ge a_n$, and
$F_b(m+2p)\ge F_b(n)$ show that the failure probability at squared
threshold
$(m-j+1)a_n^2F_b(n)/(2\beta m)$ is at most $2^{-p}$.
Use the transposed orientation instead if it supplies the larger
term in $M_n$. Since
$(m-j+1)/(2m)\ge(n-j)/(4n)\ge1/8$, we obtain
\begin{equation}\label{lin:lower-budget}
 \Pbb\{s_{\lfloor n/2\rfloor}(T+F)<h_n/(4\sqrt\beta)\}
 \le2^{-\lfloor n/8\rfloor}.
\end{equation}
The assertion is empty when $h_n=0$.
Borel--Cantelli and the finite-dilation criterion give
$\I(h)\subseteq\I(T+F)$ if $h$ has infinite support.
If $h$ has nonzero finite support, one factor has finite rank and
$T+F\ne0$ almost surely: a nonzero scalar linear test of $T$ has a
continuous density, so cannot cancel a fixed value with positive
probability. Its ideal therefore contains the finite-rank ideal
$\I(h)$. Zero factors are immediate.
Completion, Lemma~\ref{ideal:completion} with $C=0$, now proves
\eqref{ext:curved-linear-law}.

\emph{Conditional curvature.}
Fix a finite target block $U_J$ and exhaust $V$ by $V^{(N)}$.
Restricting the joint marginal potential to an observation fiber shows
that $\mu_N=\mathcal L(U_J\mid V^{(N)})$ has the same Hessian bounds.
The conditional second moments converge almost surely, giving tightness;
martingale convergence on a countable convergence-determining family
of bounded Lipschitz functions then gives
$\mu_N\Rightarrow\mu=\mathcal L(U_J\mid\mathcal F)$.
Also $m_N=\E[U_J\mid V^{(N)}]\to m=\E[U_J\mid\mathcal F]$.
Apply Lemma~\ref{curv:stability} to the translated laws
$\mu_N(\,\cdot+m_N)$, with reference covariance $I$.
Intersecting over the finite target blocks shows that
$W=U-\E[U\mid\mathcal F]$ is a curvature array on almost every
conditional fiber, with the same $\alpha,\beta$.

\emph{The innovation and its fiberwise limit.}
Set $H=KV^*D_v$. Applied to $H^*$, the first assertion gives
$H\in L^2(\mathcal S_2)$ and $\I(H)=\I(s(K)\odot v)$ almost surely.
For simple $\mathcal F$-measurable coefficients $H'$ with finite
coordinate support, conditional covariance gives
\[
 \E[\|D_uWH'\|_2^2\mid\mathcal F]
 \le\alpha^{-1}\|u\|_2^2\|H'\|_2^2.
\]
Starting with simple $\mathcal F$-measurable Hilbert--Schmidt
coefficients and using $L^2$ approximation, this estimate constructs
an $L^2(\mathcal S_2)$ limit $T$, denoted $D_uWH$.
For the finite evaluations $B_N$, let $u^{(N)},v^{(N)}$ be the
truncated weights and $\Pi_N,\Theta_N$ the noise projections. Put
\[
 H_N=\Pi_NK\Theta_NV^*D_{v^{(N)}},\qquad
 T_N=D_{u^{(N)}}WH_N.
\]
The first linear estimate gives $H_N\to H$ in $L^2(\mathcal S_2)$,
since $\Pi_NK\Theta_N\to K$ in $\mathcal S_2$ and
$v^{(N)}\to v$ in $\ell^2$. The coefficient-difference estimate gives
\begin{equation}\label{ext:innovation-convergence}
\begin{split}
 \E\|T_N-T\|_2^2
 &\le2\alpha^{-1}\|u\|_2^2\E\|H_N-H\|_2^2\\
 &\quad+2\alpha^{-1}\|u^{(N)}-u\|_2^2\E\|H\|_2^2
 \longrightarrow0.
\end{split}
\end{equation}
At each finite level the deterministic subtractions cancel, so
$T_N=B_N-\E[B_N\mid\mathcal F]$ exactly.
Since conditional expectation is an $L^2(\mathcal S_2)$ contraction,
$T=B-M$.

To identify this global limit on the conditional fibers, choose a
subsequence along which the expected squared errors in both
$H_N\to H$ and $T_N\to T$ are summable. Conditional expectation
and Tonelli give convergence in conditional $L^2$ on almost every
fiber. On such a fiber the coefficients converge, and the same
conditional covariance estimate identifies the limit with the
fixed-coefficient linear series $D_uWH$.

\emph{The conditional response.}
On each admissible fiber apply \eqref{ext:curved-linear-law}
with fixed coefficients $D_u,H$ and shift $M+F$.
Square-tail compatibility and the ideal identity for $H$ give
\eqref{ext:prediction-law}. Its event is expressed by countable
singular-value comparisons; integration against the regular
conditional kernel therefore proves the assertion for the chosen
shift $F$.
\end{proof}

\section{Exterior evaluation and covariance compensation}\label{ext:section}
\subsection{The model and its spectral shell}

A Hilbert two-form on a real separable $E$ is a skew-adjoint
Hilbert--Schmidt operator under $f\wedge g\mapsto f\otimes g-g\otimes f$;
for orthonormal increasing wedges, $\|K\|_2^2=2\|\omega_K\|_{\Lambda^2E}^2$.
Let $Q\ge0$ be trace class on $H$, and use a spectral basis of $Q$ to
write $a_n=\sqrt{q_n}$, $q_n=s_n(Q)$. For a curvature array $Z$, define
$\mathcal Zf=D_aZf$ by $L^4(H)$ limits. These exist by
\eqref{curv:convex-order}. The weight $Q$ need not be the actual output
covariance, and $\mathcal Z$ is not asserted samplewise bounded.
The full covariance of the noise is the bounded form $\Gamma$ on
$\ell^2\otimes E$ given by
\[
 \langle\Gamma(e_i\otimes f),e_j\otimes g\rangle
   =\E[Z_i(f)Z_j(g)],\qquad 0\le\Gamma\le\alpha^{-1}I.
\]
For finite-rank skew $F$, write
\[
 A_a(F)=\mathcal ZF\mathcal Z^*,\qquad
 \mathfrak c_a(F)=\E A_a(F),\qquad
 W_a(F)=A_a(F)-\mathfrak c_a(F).
\]
Here $a=(\sqrt{q_i})_i$; the notation $W_{Q,Z}$ in the introduction
denotes the same map as $W_a$.
The raw evaluation is a finite sum of vector exterior products, so
$\rank A_a(F)\le\rank F$; its mean may have infinite rank.

Fix a skew-adjoint $K\in\mathcal S_2(E)$ and choose its ordered real
spectral decomposition, completing it with zero blocks when needed:
\begin{equation}\label{ext:canonical}
 K=\sum_{\ell\ge1}K_\ell,\qquad
 K_\ell=\sigma_\ell(f_\ell\otimes g_\ell-g_\ell\otimes f_\ell),
 \qquad \kappa=s(K)=\Delta_2\sigma.
\end{equation}
Put $F_n=\sum_{\ell<n}K_\ell$, $K^{(n)}=K-F_n$, and let $P_n$ project
onto the spatial coordinates $i\ge2n-1$. Thus $F_1=0$, $P_1=I$ and
$\rank(I-P_n)\le2n-2$. We use
\begin{equation}\label{ext:profiles}
 R_n=\sum_{j\ge n}q_j,\qquad T_n=\sum_{j\ge n}\kappa_j^2,
 \qquad c_Q=(\operatorname{tr}Q)^2-\operatorname{tr}(Q^2).
\end{equation}
The two spatial factors and the noise factor combine through the
square-tail identity:
\begin{equation}\label{ext:ideal}
 \begin{gathered}
 r=a\odot\kappa\odot a,\qquad F_r(n)=R_n^2T_n,\\
 r_n^2=R_n^2T_n-R_{n+1}^2T_{n+1}
      =q_n(R_n+R_{n+1})T_n+\kappa_n^2R_{n+1}^2.
 \end{gathered}
\end{equation}
Since $\kappa_n^2\le T_n$, this gives
$\max\{\sqrt{q_nR_nT_n},\kappa_nR_n\}\le r_n
\le\sqrt3\max\{\sqrt{q_nR_nT_n},\kappa_nR_n\}$.

For finite-rank skew kernels $F$, covariance domination on linear tests gives
\begin{equation}\label{ext:trace-contraction}
 \mathfrak c_a(F)=D_aB_FD_a,\qquad
 \|B_F\|\le\alpha^{-1}\|F\|_1,\qquad
 \|P_n\mathfrak c_a(F)P_n\|
       \le\alpha^{-1}q_{2n-1}\|F\|_1.
\end{equation}
Indeed a rank-one noise kernel pairs two spatial linear tests, whose
covariance is bounded by $\alpha^{-1}$ times their norms; sum a singular
expansion of $F$. If $F$ is skew, so is $B_F$.

\subsection{Centered moments and covariance compensation}

\begin{lemma}[Polarization comparison]\label{ext:polarization}
Let $X\in\mathbb R^d$ have a centered positive density $ce^{-V}$ with
$\nabla^2V\ge\alpha I$, and let $\mathsf B$ be a symmetric bilinear map
into a finite-dimensional real normed space. For nonnegative convex $\Phi$,
\begin{equation}\label{ext:quadratic-convex}
 \E\Phi\bigl(\mathsf B(X,X)-\E\mathsf B(X,X)\bigr)
 \le\E\Phi\bigl(2\alpha^{-1}\mathsf B(g,h)\bigr),
\end{equation}
where $g,h$ are independent standard Gaussian vectors.
\end{lemma}
\begin{proof}
Take an independent copy $X'$ and put $U=(X+X')/\sqrt2$,
$V_0=(X-X')/\sqrt2$. Conditional Jensen and polarization bound the left
side by $\E\Phi(2\mathsf B(U,V_0))$. Given $U=u$, the potential of $V_0$
is $V((u+v)/\sqrt2)+V((u-v)/\sqrt2)$ plus a constant. It is even in $v$
and has Hessian at least $\alpha I$, so its mean is zero.
For fixed $u$, the map $v\mapsto\Phi(2\mathsf B(u,v))$ is convex.
Harg\'e's theorem therefore replaces the conditional variable $V_0$ by
$\alpha^{-1/2}h$, with $h$ independent of $U$. Next, independent convex-order
comparison for $X$ and $X'$ gives
\[
 \frac{X+X'}{\sqrt2}\preceq_{\rm cx}
 \frac{\alpha^{-1/2}g_1+\alpha^{-1/2}g_2}{\sqrt2}
 \stackrel{d}{=}\alpha^{-1/2}g.
\]
Indeed one replaces the first summand and then the second, keeping the
other fixed in each convex test. Finally,
$u\mapsto\E_h\Phi(2\mathsf B(u,\alpha^{-1/2}h))$ is convex, so this
second comparison gives \eqref{ext:quadratic-convex}, with independent
$g,h$. Smooth approximation gives the same statement for
$C^{1,1}_{\mathrm{loc}}$ potentials; the norm-power tests used below have
finite moments by Gaussian convex domination.
\end{proof}

\begin{lemma}[Centered evaluation and residual moments]\label{ext:centered}
The map $W_a$ extends uniquely and continuously to skew
Hilbert--Schmidt kernels, with
\[
 \E\|W_a(K)\|_2^2\le2\alpha^{-2}c_Q\|K\|_2^2.
\]
For $W_n=P_nW_a(K^{(n)})P_n$ and an absolute $C$,
\begin{equation}\label{ext:centered-residual}
 \|W_n\|_{L^n({\rm op})}\le C\alpha^{-1}r_n\qquad(n\ge1).
\end{equation}
\end{lemma}
\begin{proof}
Apply Lemma~\ref{ext:polarization} to
$\mathsf B(U,V)=\tfrac12D_d(UKV^*+VKU^*)D_d$.
For $C_d=D_dG_1KG_2^*D_d$ with independent standard Gaussian arrays,
\begin{equation}\label{ext:decoupled-upper}
 \|W_d(K)\|_{L^p({\rm op})}
 \le2\alpha^{-1}\|C_d\|_{L^p({\rm op})}\qquad(p\ge1).
\end{equation}
The Hilbert--Schmidt comparison is more precisely with
$\alpha^{-1}(C_d-C_d^*)$. For finite coefficients,
$\E|(C_d)_{ij}|^2=d_i^2d_j^2\|K\|_2^2$. If $i\ne j$, the two
entries $(C_d)_{ij}$ and $(C_d)_{ji}$ have zero covariance, while the
diagonal of $C_d-C_d^*$ vanishes. Consequently
\[
 \E\|C_d-C_d^*\|_2^2
 =2\sum_{i\ne j}d_i^2d_j^2\|K\|_2^2
 =2\left(\Bigl(\sum_i d_i^2\Bigr)^2-\sum_i d_i^4\right)\|K\|_2^2.
\]
Apply the same comparison to the difference of two skew kernels.
It gives an $L^2(\mathcal S_2)$-Cauchy sequence whenever the kernels
converge in $\mathcal S_2$. For spatial truncations, apply the entrywise
second-moment bound to the deleted entries and sum their squared weights;
that sum tends to zero because $d\in\ell^2$. Thus coordinate and spectral
approximations give the same centered extension. For a finite-rank kernel,
$L^4(H)$ continuity of $f\mapsto\mathcal Zf$ identifies this extension
with the finite sum of vector exterior products minus its expectation.

For $b=\|d\|_\infty$, $s=\|d\|_2$, $h=\|L\|_2$, $\ell=\|L\|$,
Gaussian comparison and concentration give
\begin{equation}\label{ext:five-budget}
 \|D_dG_1LG_2^*D_d\|_{L^p({\rm op})}
 \le C\bigl(bsh+s^2\ell+\sqrt p\,b^2h+
                         \sqrt p\,bs\ell+pb^2\ell\bigr).
\end{equation}
Indeed for $H_0=LG_2^*D_d$, \cite[Theorems 6.8 and 3.25]{vH} give
\[
 \|H_0\|_{L^p({\rm op})}\le C(hb+\ell s+\sqrt p\,\ell b),\qquad
 \|H_0\|_{L^p(\mathcal S_2)}\le hs+C\sqrt p\,\ell b.
\]
Condition on $G_2$ and bound the first contraction by
$C\{b\|H_0\|_{L^p(\mathcal S_2)}+(s+\sqrt p\,b)
\|H_0\|_{L^p({\rm op})}\}$. Expanding the product gives the five displayed terms.
Use $p=n$, $d=P_na$, $L=K^{(n)}$ and $nq_{2n-1}\le R_n$.
Its five terms are at most constant multiples of
$\sqrt{q_nR_nT_n}$ or $\kappa_nR_n$.
Finite compressions, $L^2$-convergent subsequences and Fatou justify
\eqref{ext:centered-residual}.
\end{proof}

A finite-rank noise head has a finite-rank raw evaluation,
but its mean may have infinite rank. To use the raw rank bound after
centering, we compensate the covariance
means on nested spatial tails. The construction makes
$P_n(\mathfrak C-\mathfrak c_a(F_n))P_n$ depend only on the kernel
blocks with index at least $n$.

\begin{lemma}[Covariance compensation]\label{ext:compensator}
Set $\gamma=\alpha^{-1}$ and $M_\ell=\mathfrak c_a(K_\ell)$.
The deterministic skew-adjoint series
\begin{equation}\label{ext:compensator-series}
 \mathfrak C=\sum_{\ell\ge1}P_{\ell+1}M_\ell P_{\ell+1}
\end{equation}
converges absolutely in $\mathcal S_2$. More precisely,
\begin{equation}\label{ext:compensator-absolute}
 \sum_{\ell\ge1}\|P_{\ell+1}M_\ell P_{\ell+1}\|_2
 \le\gamma\|K\|_2
       \Bigl(\sum_{i\ge1}(i-1)q_i^2\Bigr)^{1/2}
 \le\gamma R_1\|K\|_2.
\end{equation}
The compressed remainder satisfies
\begin{equation}\label{ext:compensator-remainder}
 \|P_n(\mathfrak C-\mathfrak c_a(F_n))P_n\|
 \le\sqrt2\gamma r_n.
\end{equation}
For fixed $a,K$ and spectral choices, the same series defines a bounded
real-linear map $\Xi\mapsto\mathfrak C_\Xi$ on bounded self-adjoint forms on
$\ell^2\otimes E$. Define its finite-rank contractions by replacing $\Gamma$ with $\Xi$
in the covariance pairing. For any two such forms,
\begin{equation}\label{ext:covariance-readout-bound}
 \|\mathfrak C_\Xi-\mathfrak C_{\widetilde\Xi}\|_2
 \le L_{a,K}\|\Xi-\widetilde\Xi\|,\qquad
 L_{a,K}=\|K\|_2\Bigl(\sum_i(i-1)q_i^2\Bigr)^{1/2}.
\end{equation}
For the actual covariance $\Gamma$, $\mathfrak C=\mathfrak C_\Gamma$.
\end{lemma}
\begin{proof}
By \eqref{ext:trace-contraction}, $M_\ell=D_aB_\ell D_a$ with
$B_\ell^*=-B_\ell$ and $\|B_\ell\|\le2\gamma\sigma_\ell$.
Schatten H\"older gives
$\|D_dBD_d\|_2\le\|D_d\|_4^2\|B\|
=(\sum_i d_i^4)^{1/2}\|B\|$.
Taking $d_i=a_i\mathbf1_{\{i\ge2\ell+1\}}$ therefore gives
\[
 \|P_{\ell+1}M_\ell P_{\ell+1}\|_2
 \le2\gamma\sigma_\ell
       \Bigl(\sum_{i\ge2\ell+1}q_i^2\Bigr)^{1/2}.
\]
Cauchy--Schwarz, $2\sum_\ell\sigma_\ell^2=\|K\|_2^2$ and
$\sum_\ell\sum_{i\ge2\ell+1}q_i^2
\le\tfrac12\sum_i(i-1)q_i^2$ prove the first bound.
Monotonicity gives
$\sum_i(i-1)q_i^2\le\sum_{j<i}q_iq_j\le R_1^2$.
For $\ell<n$, nested projections give
$P_nP_{\ell+1}=P_n$; for $\ell\ge n$, they give
$P_nP_{\ell+1}=P_{\ell+1}$. Therefore
\[
 P_n(\mathfrak C-\mathfrak c_a(F_n))P_n
 =\sum_{\ell\ge n}P_{\ell+1}M_\ell P_{\ell+1}.
\]
Its norm is at most
\[
 2\gamma\sum_{\ell\ge n}\sigma_\ell q_{2\ell+1}
 \le\sqrt2\gamma\sqrt{T_{2n-1}q_{2n+1}R_{2n+1}}
 \le\sqrt2\gamma\sqrt{T_nq_nR_n}\le\sqrt2\gamma r_n.
\]
All sums are justified by absolute convergence. The contraction estimate
uses only $\|\Xi\|$, not positivity, so all these bounds hold with
$\gamma=\|\Xi\|$. Applying the absolute-sum bound to
$\Xi-\widetilde\Xi$ proves \eqref{ext:covariance-readout-bound} and linearity.\qedhere
\end{proof}

For the spectral tail $Q_{\ge n}$ with eigenvalues
$q_n,q_{n+1},\ldots$, one has
$E_n=\operatorname{tr}(\Lambda^2Q_{\ge n})$, where
$\Lambda^2Q_{\ge n}(u\wedge v)=Q_{\ge n}u\wedge Q_{\ge n}v$.
The next profile decomposes the product of this exterior trace
and the noise tail $T_n$ into shells.

\begin{proposition}[Exterior-square generator of the fluctuation ideal]
\label{ext:exterior-shell}
Let $q=(q_n)$ and $\kappa=(\kappa_n)$ be nonnegative, nonincreasing,
with $\sum_nq_n<\infty$ and $\sum_n\kappa_n^2<\infty$.
Put $R_n=\sum_{j\ge n}q_j$, $T_n=\sum_{j\ge n}\kappa_j^2$ and
$E_n=\sum_{n\le i<j}q_iq_j$; define
\[
 \begin{aligned}
 r_n^2&=R_n^2T_n-R_{n+1}^2T_{n+1},\\
 e_n^2&=E_nT_n-E_{n+1}T_{n+1}
       =\sum_{\substack{i,j,k\ge1,\ i<j\\
                         \min\{i,j,k\}=n}}q_iq_j\kappa_k^2.
 \end{aligned}
\]

Then $e,r$ are nonnegative, nonincreasing null sequences, and
\begin{equation}\label{ext:exterior-one-step}
 e_n^2\le\tfrac12r_n^2\quad(n\ge1),\qquad
 r_n^2\le2e_{n-1}^2\quad(n\ge2).
\end{equation}
If $q_2>0$ and $\kappa_1>0$, then $\I(e)=\I(r)$.
If $q_2=0$, then $e\equiv0$ (and the spatial alternating
evaluation vanishes).
\end{proposition}
\begin{proof}
Let $S=R_{n+1}$. Expanding the differences yields
\[
 r_n^2=(q_n^2+2q_nS)T_n+S^2\kappa_n^2,
 \qquad e_n^2=q_nS\,T_n+E_{n+1}\kappa_n^2.
\]
Because $2E_{n+1}\le S^2$, the first inequality follows.
Since $q$ is nonincreasing,
$S^2\le2E_{n+1}+q_nS\le2E_n$ and
$(q_n^2+2q_nS)T_n\le2q_{n-1}R_nT_{n-1}$ for $n\ge2$.
Using $\kappa_n\le\kappa_{n-1}$ gives
$r_n^2\le2e_{n-1}^2$. Both profiles decrease to zero
by the expansions and summability assumptions.
If $q_2\kappa_1>0$, then $e_1>0$, and the one-step comparison
implies $r_n\le C e_{\lceil n/2\rceil}$ for a constant $C<\infty$.
The standard finite-dilation criterion gives $\I(r)\subseteq\I(e)$;
the first comparison gives the reverse inclusion.
\end{proof}

\subsection{Complete response}

\begin{lemma}[A reducing decomposition preserving the ideal]\label{ext:splitting}
For a positive compact $Q$ of rank at least two, there is an orthogonal
projection $P$ commuting with $Q$ such that both restrictions $Q_+,Q_-$ to
$PH,(I-P)H$ are nonzero and
\[
 \I(Q_+^{1/2})=\I(Q^{1/2})=\I(Q_-^{1/2}).
\]
\end{lemma}
\begin{proof}
Divide an ordered positive spectral basis alternately between the two
subspaces, splitting eigenspaces when necessary. If $a=s(Q^{1/2})$, their
spectra are $a_n^+=a_{2n-1}$ and $a_n^-=a_{2n}$, with zeros appended, and
\begin{equation}\label{ext:subsequence}
 a^-\le a^+\le a,\qquad
 a_n\le\max\{1,a_1/a_2\}\,a^-_{\lceil n/3\rceil}.
\end{equation}
For $n\ge2$ use $2\lceil n/3\rceil\le n$; the constant handles $n=1$.
The ideal criterion proves the assertion, including finite rank.
\end{proof}

\begin{maintheorem}[Exterior geometry and intrinsic affine response]
\label{ext:classification}
Assume \eqref{ext:array-curvature}, $\rank Q\ge2$ and $K\ne0$.
The centered alternating evaluation has a complete deterministic
response with unique intrinsic pair
$(J,C+J)=(\I(e),\mathfrak C+\I(e))$.
For every fixed compact operator $B$,
\begin{equation}\label{ext:classification-law}
 \I(W_a(K)+B)=\I(e)+\I(B-\mathfrak C)
 \quad\text{almost surely}.
\end{equation}
In particular $Y_0=W_a(K)+\mathfrak C$ has
\begin{equation}\label{ext:compensated-budget}
 \|s_{6n-5}(Y_0)\|_{L^n}\le C\alpha^{-1}r_n.
\end{equation}
\end{maintheorem}
\begin{proof}
\emph{Lower saturation.}
For every fixed compact $B$, we prove
\begin{equation}\label{ext:universal-lower}
 \I(r)\subseteq\I(W_a(K)+B)\qquad\text{almost surely}.
\end{equation}
Use the odd/even split of Lemma~\ref{ext:splitting}, with
$u_n=a_{2n-1}$, $v_n=a_{2n}$, and let $V$ be that block of $W_a(K)$.
Its centered coordinate evaluations converge in $L^2(\mathcal S_2)$,
with deterministic subtractions. The prediction law
\eqref{ext:prediction-law}, applied with the fixed shift $B_{+-}$,
gives $\I(u\odot\kappa\odot v)\subseteq\I(V+B_{+-})$
almost surely, even when the raw means diverge.
Square-tail compatibility and \eqref{ext:subsequence} identify the left
side as $\I(r)$; compression proves \eqref{ext:universal-lower}.

\emph{Compensated upper membership.}
Recall $Y_0=W_a(K)+\mathfrak C$ and set
\[
 R_n^{\rm fin}=(I-P_n)Y_0+P_nY_0(I-P_n)+P_nA_a(F_n)P_n.
\]
Each of the first two terms has rank at most $2n-2$, as does the
last term by the raw rank bound. Thus $\rank R_n^{\rm fin}\le6n-6$.
Using the raw head here is essential: its covariance mean can have
infinite rank. Linearity of the centered evaluation gives, on one
event for all $n$,
\begin{equation}\label{ext:compensated-decomposition}
 Y_0=R_n^{\rm fin}+W_n+
       P_n(\mathfrak C-\mathfrak c_a(F_n))P_n.
\end{equation}
The last two terms have $L^n$ operator norm at most
$C\alpha^{-1}r_n$, by Lemmas~\ref{ext:centered} and
\ref{ext:compensator}. The rank bound therefore gives
\eqref{ext:compensated-budget}.
Markov at $eC\alpha^{-1}r_n$ and Borel--Cantelli imply
$Y_0\in\I(r)$; zero tails give the corresponding finite-rank bounds
directly.

\emph{Affine completion.}
We have proved exactly the two inputs of Lemma~\ref{ideal:completion}:
$W_a(K)+\mathfrak C\in\I(r)$ and
$\I(r)\subseteq\I(W_a(K)+B)$ for each fixed $B$.
The lemma gives the complete response, and
Proposition~\ref{ext:exterior-shell} identifies $\I(r)=\I(e)$.
Theorem~\ref{ideal:response-structure} identifies the pair uniquely.
The moment estimate remains in the profile $r$ for which it was
proved; only ideal equivalence is claimed for $e$.
\end{proof}

\subsection{Raw limits and equivalence of responses}
For a fixed noise basis, let $E_m$ and $\Pi_n$ be the first-$m$
spatial and first-$n$ noise projections. Along any prescribed
exhaustion $m_N,n_N\uparrow\infty$, define
\begin{equation}\label{ext:raw-cutoff}
 A_N=E_{m_N}D_aZ\Pi_{n_N}K\Pi_{n_N}Z^*D_aE_{m_N}.
\end{equation}
\begin{corollary}[Raw limits and their response]
\label{ext:raw-classification}
Under Theorem~\ref{ext:classification}, $A_N$ converges in
$L^2(\mathcal S_2)$ if and only if its deterministic means converge
in $\mathcal S_2$. If $\E A_N\to M$, its limit is
$A=W_a(K)+M$, and for each fixed compact $D$,
\begin{equation}\label{ext:raw-law}
 \I(A+D)=\I(e)+\I(D+M-\mathfrak C)\quad\text{almost surely}.
\end{equation}
In particular the raw response is $\I(e)+\I(D)$ exactly when
$M-\mathfrak C\in\I(e)$.
\end{corollary}
\begin{proof}
The centered finite section is
$E_{m_N}W_a(\Pi_{n_N}K\Pi_{n_N})E_{m_N}$.
The continuity in Lemma~\ref{ext:centered} and dominated convergence
of spatial compressions show that it converges in
$L^2(\mathcal S_2)$ to $W_a(K)$. Expectation is continuous in this
space, proving the existence assertion. Substitute $B=M+D$ in
\eqref{ext:classification-law}. The last assertion follows by taking
$D=\mathfrak C-M$ in the claimed additive response.
\end{proof}

Even two spatial directions can couple the kernel to a bounded skew
covariance block. The universal raw domain is therefore governed by
trace-class/bounded-operator duality.

\begin{proposition}[The universal domain of raw evaluation]
\label{ext:universal-raw-domain}
Fix a positive trace-class $Q$ of rank at least two and a skew
Hilbert--Schmidt kernel $K$.
Along any prescribed exhaustion in \eqref{ext:raw-cutoff}, raw
evaluation converges in $L^2(\mathcal S_2)$ for every curvature
array if and only if $K\in\mathcal S_1$. In that case the
convergence holds in $L^2(\mathcal S_1)$.
If $K\notin\mathcal S_1$, failure already occurs for centered
Gaussian arrays with covariance arbitrarily close to the identity.
\end{proposition}
\begin{proof}
For a rank-one noise kernel $f\otimes g$, convex comparison gives
\[
 \|(\mathcal Zf)\otimes(\mathcal Zg)\|_{L^2(\mathcal S_1)}
 \le \sqrt3\,\alpha^{-1}(\operatorname{tr}Q)\|f\|\|g\|.
\]
Summing a singular expansion gives a continuous linear raw map
$\mathcal S_1(E)\to L^2(\mathcal S_1(H))$.
The noise compressions converge in trace norm to $K$, and spatial
compressions of their limiting raw evaluation converge in
$L^2(\mathcal S_1)$ by dominated convergence. This proves
sufficiency along every prescribed exhaustion.

Conversely, put $K_N=\Pi_{n_N}K\Pi_{n_N}$.
If their trace norms were bounded, singular-value convergence and
Fatou would imply $K\in\mathcal S_1$. Thus
$\sup_N\|K_N\|_1=\infty$.
On the real Banach space of bounded skew-adjoint operators on $E$,
the functionals $L\mapsto\operatorname{tr}(K_N^*L)$ have norm
$\|K_N\|_1$: the skew polar partial isometry attains it.
The uniform boundedness principle supplies a bounded skew $L$
with $\sup_N|\operatorname{tr}(K_N^*L)|=\infty$.
Scaling gives $\|L\|<\varepsilon$ for any prescribed
$0<\varepsilon<1$, without removing unboundedness.
Give the first two Gaussian rows joint covariance
\[
 \begin{pmatrix}I&L\\-L&I\end{pmatrix}
\]
and take the other rows independently standard.
The full covariance lies between $(1-\varepsilon)I$ and
$(1+\varepsilon)I$, so its finite marginals satisfy the curvature
hypothesis with constants $(1+\varepsilon)^{-1}$ and
$(1-\varepsilon)^{-1}$. Once both spatial coordinates are retained,
the first off-diagonal raw mean is
$\sqrt{q_1q_2}\operatorname{tr}(K_N^*L)$.
It is unbounded; continuity of expectation precludes
$L^2(\mathcal S_2)$ convergence. The argument uses the prescribed
noise exhaustion throughout.
\end{proof}

\begin{corollary}[Exactly which covariance information is retained]
\label{ext:response-equivalence}
For fixed $Q,K$ as in Theorem~\ref{ext:classification}, let two
curvature arrays have covariances $\Gamma_1,\Gamma_2$ and centered
evaluations $W_1,W_2$. Write $J=\I(e)$ and
$C_i=\mathfrak C_{\Gamma_i}$. Their complete responses coincide
if and only if
\begin{equation}\label{ext:covariance-class}
 C_1-C_2\in J.
\end{equation}
On every prescribed coupling this is also equivalent to
$W_1-W_2\in J$ almost surely. The coset $C_i+J$ does not depend
on the allowed spectral choices in its construction. Equal array
covariances are sufficient for equal responses and for the same
raw-convergence behavior along a fixed exhaustion.
\end{corollary}
\begin{proof}
Uniqueness in Theorem~\ref{ideal:response-structure} proves the
response criterion; necessity can already be read at $B=C_1$.
Subtracting $W_i+C_i\in J$ proves both directions of the coupling
criterion. Different spectral choices classify the same response,
so uniqueness gives the same coset. The compensator and every
finite raw mean are covariance contractions, giving the final
assertion. If $\rank Q\le1$ or $K=0$, alternation instead gives
$W=A=0$ and the response is simply $\I(B)$.
\end{proof}

\subsection{The interval of attainable raw ideals}\label{ext:interval-section}

For a fixed law with $K\in\mathcal S_1$, \eqref{ext:raw-law} identifies
the minimizing coset as $-\Delta+J$. At zero shift the ideal is
$J+\I(\Delta)$. Which values can it take as the joint law of the rows varies?

\begin{maintheorem}[The attainable interval for a fixed pair $Q,K$]
\label{ext:raw-interval}
Let $Q\ge0$ be trace class of rank at least two, and let $K$ be a
nonzero skew-adjoint trace-class kernel. Put $J=\I(e)$ and
$\Delta=\E A_Z-\mathfrak C_\Gamma$. With \eqref{ext:profiles} write
\begin{equation}\label{ext:interval-endpoint}
 S_n=\sum_{j\ge n}\kappa_j,\qquad
 p_n=R_nS_n-R_{n+1}S_{n+1}=q_nS_n+\kappa_nR_{n+1}.
\end{equation}
For every centered curvature array $Z$, its raw evaluation $A_Z$
satisfies
\begin{equation}\label{ext:interval-bound}
 s_{4n-3}(\Delta)\le 2\alpha^{-1}q_{2n-1}S_{2n-1},
 \qquad J\subseteq\I(A_Z)\subseteq\I(p)\quad\text{almost surely}.
\end{equation}
Allowing the curvature constants to depend on the law, the attainable
ideals are precisely the principal ideals in this interval. Take any
nonincreasing null sequence $u\ge0$ with
$J\subseteq\I(u)\subseteq\I(p)$, and any $0<\varepsilon<1$.
For these same $Q,K$ there is a centered Gaussian array with standard
Gaussian row marginals and $\|\Gamma-I\|\le\varepsilon$ for which
$\I(A_Z)=\I(u)$ almost surely. In the chosen spectral bases put
\[
 J_0=\begin{pmatrix}0&1\\-1&0\end{pmatrix},\qquad
 \delta_i=a_{2i-1}a_{2i}S_{2i-1},\qquad
 M_u=\varepsilon\bigoplus_{i\ge1}\min(u_i,\delta_i)J_0 .
\]
This array has $\mathfrak C_\Gamma=0$ and $\E A_Z=M_u$. For each fixed
compact $D$,
\begin{equation}\label{ext:interval-shifts}
 \I(A_Z+D)=J+\I(D+M_u)\quad\text{almost surely}.
\end{equation}
The shifts attaining $J$ almost surely form the coset $-M_u+J$.
\end{maintheorem}

\begin{proof}
By Proposition~\ref{ext:universal-raw-domain} the raw map is continuous
from trace class into $L^2(\mathcal S_1)$. In particular $M=\E A_Z=\sum_\ell M_\ell$.
Write $\gamma=\alpha^{-1}$. By \eqref{ext:trace-contraction},
$\|M_\ell\|_1\le2\gamma\sigma_\ell R_1$, so both $M$ and
\eqref{ext:compensator-series} converge absolutely in trace norm.
On the spatial tail $P_n$, their first $n-1$ terms cancel:
\[
 P_n\Delta P_n=\sum_{\ell\ge n}
 (P_nM_\ell P_n-P_{\ell+1}M_\ell P_{\ell+1}).
\]
Each summand has norm at most $4\gamma\sigma_\ell q_{2n-1}$.
Since $2\sum_{\ell\ge n}\sigma_\ell=S_{2n-1}$ and
$\rank(\Delta-P_n\Delta P_n)\le 2\rank(I-P_n)\le4n-4$,
\eqref{var:approximation} proves the singular-value bound.
Finite dilation now gives $\I(\Delta)\subseteq\I((q_nS_n)_n)$.

The nonnegative sequence $p$ is nonincreasing by \eqref{ext:interval-endpoint}
and sums to $R_1S_1$. From $T_n\le\kappa_n S_n\le S_n^2$ and
\eqref{ext:ideal} we get $\kappa_nR_{n+1}\le r_n\le p_n$ by squaring.
Thus $q_nS_n\le p_n\le q_nS_n+r_n$. Since $J=\I(r)$ by
Proposition~\ref{ext:exterior-shell}, we have
\[
 \I(p)=J+\I((q_nS_n)_n).
\]
Together with \eqref{ext:raw-law} this proves both ideal inclusions.

To realize $\I(u)$, put $t_i=\min(u_i,\delta_i)$.
We have $\delta_n\le q_nS_n$, and
$q_nS_n\le\max(1,a_1/a_2)\delta_{\lceil n/3\rceil}$, since
$2\lceil n/3\rceil\le n$ for $n\ge2$. Hence
$\I(\delta)=\I((q_nS_n)_n)$. Coordinatewise minimum distributes
over maximum, so Lemma~\ref{var:idealcalculus} gives
\[
 \I(u)=\I(u)\cap(J+\I(\delta))=J+(\I(u)\cap\I(\delta))=J+\I(t).
\]
This uses only ideal inclusion for $u$.

Using zero extension when needed, decompose the spatial and noise coordinates into the pairs used in
\eqref{ext:canonical} and in the definition of $P_n$. On spatial
pair $i$ tensored with noise pair $\ell$, prescribe the covariance
\[
 \Gamma_{i\ell}=I_4+c_{i\ell}J_0\otimes J_0,\qquad
 c_{i\ell}=\mathbf1_{\{\ell\ge i\}}\varepsilon\,t_i/\delta_i,
\]
with zero quotient when $\delta_i=0$. Use independent blocks and
independent standard noise on the orthogonal complement.
The matrix $J_0\otimes J_0$ is a symmetric involution with zero
diagonal spatial blocks. The spectrum of $\Gamma$ therefore lies
in $[1-\varepsilon,1+\varepsilon]$, and its compression to each
single row is the identity. This covariance defines a Gaussian field
in any prescribed noise basis, with curvature bounds
$(1+\varepsilon)^{-1}$ and $(1-\varepsilon)^{-1}$.

Contracting $K_\ell=\sigma_\ell J_0$ with the off-diagonal row
covariance $c_{i\ell}J_0$ gives $2\sigma_\ell c_{i\ell}$.
Thus the $i$th spatial block of $M_\ell$ is
$2a_{2i-1}a_{2i}\sigma_\ell c_{i\ell}J_0$. Only $i\le\ell$
can contribute, so $P_{\ell+1}M_\ell P_{\ell+1}=0$ for every $\ell$.
Summing in $\ell$ gives $M_u$ and the compensator vanishes.
Each $t_i$ occurs twice in $s(M_u)/\varepsilon$, so $\I(M_u)=\I(t)$.
Apply \eqref{ext:raw-law} at $D=0$ and then at every fixed $D$ to
obtain the realization and \eqref{ext:interval-shifts}.
\end{proof}

For example, take $q_n=n^{-2}$ and
$\kappa_{2i-1}=\kappa_{2i}=1/[i\log^2(ei)]$.
Integral comparison gives $R_n\asymp n^{-1}$,
$S_n\asymp 1/\log(en)$ and $T_n\asymp 1/[n\log^4(en)]$.
Substitution in \eqref{ext:ideal} and \eqref{ext:interval-endpoint}
yields $r_n\asymp n^{-2}\log^{-2}(en)$ and
$p_n\asymp n^{-2}\log^{-1}(en)$.
Thus $u^{(\theta)}_n=n^{-2}\log^{\theta-2}(en)$, $0\le\theta\le1$,
gives a strictly increasing chain of attainable ideals for one pair
$Q,K$, since fixed index dilations preserve the logarithmic exponent.
All these profiles belong to $\ell^s$ exactly for $s>1/2$. At $s=1/2$
the series is $\sum_n n^{-1}\log^{-1+\theta/2}(en)$ and diverges.

\subsection{Gaussian observations and the minimal residual ideal}

For independent standard Gaussian rows, write $A_Q(K)=W_a(K)$;
the raw means vanish, so this is also the raw Hilbert--Schmidt
evaluation. Fix a deterministic orthogonal projection $P$ on $E$ and set
\[
 \mathcal G_P=\sigma\{\mathcal Zf:f\in\operatorname{ran}P\},
 \qquad M_P=\E[A_Q(K)\mid\mathcal G_P].
\]

\begin{theorem}[Complete response under Gaussian observations]
\label{ext:gaussian-information}
Let $\rank Q\ge2$ and $K$ be skew Hilbert--Schmidt. Put
$L=(I-P)K$ and $J_P=\I(a\odot s(L)\odot a)$.
Then $M_P=A_Q(PKP)$, and every chosen
$\mathcal G_P$-measurable compact random operator $D$ satisfies
\begin{equation}\label{ext:observed-response}
 \I(A_Q(K)+D)=J_P+\I(M_P+D)\quad\text{almost surely}.
\end{equation}
The least attainable ideal is $J_P$; its minimizing corrections are
exactly those satisfying $D+M_P\in J_P$ almost surely.
No finite-dimensionality or reducing assumption on $P$ is required;
the probability-one event is for each chosen correction separately.
\end{theorem}
\begin{proof}
Put $R=I-P$ and $K_0=K-PKP$. Orthogonal Gaussian decomposition
gives the conditional-mean identity. Moreover
\[
 K_0=L-PL^*,\qquad L=RK_0,
\]
so $\I(K_0)=\I(L)$. Square-tail compatibility and
Theorem~\ref{ext:classification} (or the trivial zero-kernel case)
give $A_Q(K)-M_P=A_Q(K_0)\in J_P$ almost surely.
If $L=0$, this proves the result.

For the lower inclusion use the spatial split of
Lemma~\ref{ext:splitting}, with $u_n=a_{2n-1}$ and $v_n=a_{2n}$.
Let $V$ be the corresponding off-diagonal block of $A_Q(K)$ and
condition on $\mathcal F=\mathcal G_P\vee\sigma(Z_-)$, where
$Z_-$ denotes all negative spatial rows. The positive residual rows
$Y_+=Z_+|_{RE}$ are independent standard Gaussian rows, and
\[
 V-\E[V\mid\mathcal F]=D_uY_+H,
 \qquad H=LZ_-^*D_v.
\]
These are Hilbert--Schmidt limits, with
$\E\|H\|_2^2=\|L\|_2^2\|v\|_2^2$.
Finite-coordinate approximation in a basis adapted to $P$ and the
conditional subsequence
argument in Lemma~\ref{ext:curved-linear} identify the innovation
on almost every fiber. The unconditional linear law gives
$\I(H)=\I(s(L)\odot v)$ almost surely. Applying the same law on
each fiber with shift $\E[V\mid\mathcal F]+D_{+-}$ gives
\[
 \I(u\odot s(L)\odot v)
 \subseteq\I(V+D_{+-})\subseteq\I(A_Q(K)+D).
\]
The left side equals $J_P$ by the splitting lemma and square-tail
compatibility. Borel ideal comparison permits conditioning for each
chosen $D$. Subtraction, using $A_Q(K)-M_P\in J_P$, now proves
\eqref{ext:observed-response}.
\end{proof}

For every algebraic two-sided ideal $N$ and every
$\mathcal G_P$-measurable compact $B$, the response formula gives
\[
 A_Q(K)-B\in N\ \text{a.s.}
 \quad\Longleftrightarrow\quad
 J_P\subseteq N\ \text{and}\ B-M_P\in N\ \text{a.s.}
\]
Membership is asserted on measurable events of probability one,
as in Section~\ref{ideal:section}; take $D=-B$ in
\eqref{ext:observed-response}. Thus $M_P+J_P$ is the least affine
carrier with an observation-measurable center and a deterministic
ideal direction, the observed counterpart of $-C+J$ in
Theorem~\ref{ideal:response-structure}. No measurable quotient is needed.
Exact recovery occurs precisely when $(I-P)K=0$.
For a fixed pair $P\le P'$, apply the preceding criterion at $P'$
with the $\mathcal G_{P'}$-measurable center $B=M_P$ and $N=J_P$.
It gives $M_{P'}+J_{P'}\subseteq M_P+J_P$ almost surely, hence
$J_{P'}\subseteq J_P$.

Write $J_0=\I(a\odot s(K)\odot a)$. The case $0\le P$ also gives
$M_P\in J_0$ almost surely. Consequently, the unobserved ideal
identity $\I(A_Q(K)+D)=J_0+\I(D)$ remains valid almost surely for
each chosen $\mathcal G_P$-measurable compact correction $D$
if and only if $J_P=J_0$: sufficiency follows from
\eqref{ext:observed-response}, and necessity follows by taking $D=-M_P$.
Finite-rank $P$ and infinite-rank $K$ satisfy this criterion because
$\I((I-P)K)=\I(K)$.
The one-sided kernel is essential. If $P$ observes every $f_j$ but
none of the $g_j$ in \eqref{ext:canonical}, then
$PKP=(I-P)K(I-P)=0$, whereas
$s((I-P)K)=(\sigma_j)_j$ and $s(K)=\Delta_2\sigma$ generate the same
ideal. Such observations leave the response ideal unchanged.

In contrast, suppose $\sigma_j\asymp j^{-1}$ and $q_n\asymp n^{-2}$,
and observe the complete canonical pairs whose indices are not squares.
Then $\lambda_n:=s((I-P)K)_n\asymp n^{-2}$ and $R_n\asymp n^{-1}$.
Writing $T_n^P=\sum_{k\ge n}\lambda_k^2$, the original and remaining
kernel tails satisfy $T_n\asymp n^{-1}$ and $T_n^P\asymp n^{-3}$.
The two contributions in \eqref{ext:ideal} consequently have orders
$n^{-4}$ for the original squared shell and $n^{-6}$ for the remaining
one. Thus
\[
 \I(A_Q(K))=\I((n^{-2})_n),\qquad
 J_P=\I((n^{-3})_n)\subsetneq\I((n^{-2})_n).
\]
With the same observation scheme but $q_n=e^{-n}$,
\eqref{ext:ideal} gives original and residual shell profiles comparable
to $e^{-n}n^{-1/2}$ and $e^{-n}n^{-3/2}$, respectively.
Both generate $\I((e^{-n})_n)$, since a fixed index dilation absorbs
these polynomial factors.
Thus the complete ideal response is unchanged although
$\I(s((I-P)K))\subsetneq\I(s(K))$, even with infinite-rank $Q$.

\section{Wiener area and covariance tails}
\label{area:section}

We now compute the fixed-shift response of Hilbert-space Wiener area.
Exterior-dilation estimates and finite-dimensional singular-value
distributions are studied in \cite{ZhaoExterior} and \cite{DCO},
respectively.

Let $W=(W_t)_{0\le t\le1}$ be a Wiener process in a real separable
Hilbert space with trace-class covariance $Q\ge0$.  Define its
Hilbert--Schmidt-valued skew area on $[0,1]$ by
\begin{equation}\label{area:def}
 A_Q=\frac12\int_0^1
 \bigl(W_t\otimes dW_t-dW_t\otimes W_t\bigr).
\end{equation}
Write $q_n=s_n(Q)$, $R_n=\sum_{k\ge n}q_k$ and
\begin{equation}\label{area:b}
 b_n=\frac{R_n}{n},\qquad n\ge1.
\end{equation}
The trace-class assumption makes the It\^o integral in
\eqref{area:def} well defined in $\mathcal S_2$~\cite{DPZ2014}.

Represent $W$ using independent isonormal rows
$Z_i(f)=\int_0^1f(t)\,d\beta_i(t)$ on $L^2(0,1)$.  The
antisymmetric time kernel is the skew Hilbert--Schmidt operator
\[
 (Kf)(s)=\frac12\int_0^1\operatorname{sgn}(t-s)f(t)\,dt.
\]
Riemann--It\^o sums and the It\^o isometry identify the centered
alternating evaluation $W_{Q,Z}(K)$ with $A_Q$.
On complexification, for $\nu=(2j-1)\pi$ one has
$K(e^{i\nu\,\cdot})(s)=i e^{i\nu s}/\nu$.
The positive and negative odd Fourier frequencies form an orthonormal
basis, giving $\kappa_{2j-1}=\kappa_{2j}=1/[\pi(2j-1)]$.

\begin{lemma}[The Brownian time shell is a tail average]
\label{area:shell-lemma}
Let $\kappa_n=[\pi(2\lceil n/2\rceil-1)]^{-1}$,
$T_n=\sum_{k\ge n}\kappa_k^2$, and put
$r_n^2=R_n^2T_n-R_{n+1}^2T_{n+1}$.
For every decreasing nonnegative summable $q$ with $q_2>0$,
\begin{equation}\label{area:shell-ideal}
 \I(r)=\I(b).
\end{equation}
Moreover, for every such $q$,
\begin{equation}\label{area:regularity-equivalence}
 \I(q)=\I(b)
 \quad\Longleftrightarrow\quad
 \sup_{n:q_n>0}\frac{R_n}{nq_n}<\infty.
\end{equation}
For finite support the ratio condition is interpreted together with
$R_n=0$ whenever $q_n=0$; it automatically holds.
\end{lemma}
\begin{proof}
The explicit $\kappa$ has, with absolute positive constants,
$\kappa_n\asymp n^{-1}$ and $T_n\asymp n^{-1}$.
The shell expansion is
\[
 r_n^2=q_n(R_n+R_{n+1})T_n+\kappa_n^2R_{n+1}^2.
\]
Consequently
\[
 c\left(\frac{q_nR_n}{n}
       +\frac{R_{n+1}^2}{n^2}\right)
 \le r_n^2
 \le C\left(\frac{q_nR_n}{n}
       +\frac{R_n^2}{n^2}\right).
\]
For $R_n>0$, set $x=R_{n+1}/R_n\in[0,1]$.  Since
$(1-x)+x^2\ge3/4$, the lower inequality gives
$r_n\ge c_0R_n/n$.  If $R_n=0$ both sides vanish.
For $n\ge2$, set $m=\lceil n/2\rceil$.
Monotonicity of $q$ gives
$R_m\ge(n-m+1)q_n\ge(n/2)q_n$.
Thus the upper inequality yields $r_n\le C_0R_m/m$.
The first index is covered by changing $C_0$, since $q_2>0$.
The finite-dilation criterion proves \eqref{area:shell-ideal}.

The same monotonicity bound shows $q_n\le b_{\lceil n/2\rceil}$,
so $\I(q)\subseteq\I(b)$.  If $R_n\le Cnq_n$, the reverse
inclusion is immediate.  Conversely, if $\I(b)\subseteq\I(q)$,
then $b_{Mk}\le Cq_k$ for some fixed $M\ge1$ and $C$.
Decomposing the tail at $Mk$ gives
\[
 R_k=\sum_{j=k}^{Mk-1}q_j+R_{Mk}
 \le(M-1+CM)kq_k.\qedhere
\]
\end{proof}

\begin{corollary}[Complete Wiener-area response and exact covariance criterion]
\label{area:main}
If $\rank Q\ge2$, then for every fixed deterministic compact
operator $D$, not necessarily skew-adjoint,
\begin{equation}\label{area:response}
 \I(A_Q+D)=\I\bigl((R_n/n)_{n\ge1}\bigr)+\I(D)
 \quad\text{almost surely}.
\end{equation}
Its intrinsic affine-response pair is
$\bigl(\I(b),0+\I(b)\bigr)$, and in particular
\begin{equation}\label{area:comparison}
 \I(Q)\subseteq\I(A_Q),\qquad
 \I(Q)=\I(A_Q)
 \ \Longleftrightarrow\ 
 R_n\le Cnq_n\ \text{for all }n\ge1
 \text{ and some }C<\infty.
\end{equation}
If $\rank Q\le1$, then $A_Q=0$ almost surely and
$\I(A_Q+D)=\I(D)$.
The same identities hold for every fixed nondegenerate interval
$[s,t]$ after recentering the path at $s$.
When $\rank Q\ge2$, \eqref{area:response} also holds for every
chosen compact random $D$ measurable with respect to
\begin{equation}\label{area:finite-information}
 \mathcal G_P=\sigma\left\{\int_0^1 f(t)\,dW_t:
                              f\in\operatorname{ran}P\right\},
\end{equation}
where $P$ is any fixed finite-rank orthogonal projection on $L^2(0,1)$.
Thus finitely many Hilbert-valued temporal observations, including
values at finitely many times, permit no improvement of $\I(b)$ by
any measurable compact correction. The event is for each chosen
correction separately.

\end{corollary}
\begin{proof}
Independence of the rows makes the covariance compensator zero.
Theorem~\ref{ext:classification} and
Proposition~\ref{ext:exterior-shell} therefore give
$\I(A_Q+D)=\I(r)+\I(D)$ for each fixed $D$.
Lemma~\ref{area:shell-lemma} proves \eqref{area:response}
and \eqref{area:comparison}.  Rank at most one forces the
alternating spatial area to vanish.  Brownian scaling multiplies
any fixed-interval area by a nonzero deterministic scalar in law,
which does not change a principal ideal or its full fixed-shift response.

Since the Brownian kernel has infinite rank, the finite-rank observation
case of Theorem~\ref{ext:gaussian-information} gives
\eqref{area:response} for every chosen
$\mathcal G_P$-measurable compact correction $D$.

\end{proof}

\begin{corollary}[Nuclearity of the area ideal]\label{area:nuclearity}
Let $Q\ge0$ be trace class of infinite rank, with $q=s(Q)$ and
$b_n=R_n/n$. Then $\I(q)$ satisfies \eqref{fs:criterion} if and
only if $\I(b)$ does. In that case $\I(A_Q)=\I(Q)$ almost surely.
Thus the scalar-generation and faithful-geometry alternatives for
the area ideal are determined by the same nuclearity test on $q$.
\end{corollary}
\begin{proof}
The dyadic criterion gives $a_{Mn}/a_n\to0$ for some fixed $M$
whenever $a$ satisfies \eqref{fs:criterion}.
For $a=q$, choose $\rho<1/M$ with $q_{Mn}\le\rho q_n$
eventually. Summing the blocks $[M^kn,M^{k+1}n)$ gives
$R_n\le(M-1)nq_n/(1-M\rho)$.
For $a=b$, choose $b_{Mn}\le b_n/(2M)$ eventually. Then
$R_{Mn}\le R_n/2$, so
$R_n\le2\sum_{j=n}^{Mn-1}q_j\le2(M-1)nq_n$.
Either hypothesis gives $b_n\le Cq_n$; the reverse ideal inclusion
always follows from $q_n\le b_{\lceil n/2\rceil}$.
Hence $\I(q)=\I(b)$, and invariance of \eqref{fs:criterion}
proves the equivalence. Apply \eqref{area:response} and
Theorem~\ref{fs:intrinsic-dichotomy}.
\end{proof}

\begin{example}[A logarithmic separation]
\label{area:logarithmic-example}
Let $q_n\asymp[n(\log(en))^\gamma]^{-1}$ with $\gamma>2$.
Then
\[
 R_n\asymp(\log(en))^{1-\gamma},\qquad
 b_n\asymp[n(\log(en))^{\gamma-1}]^{-1}.
\]
For every fixed $m$, the ratio $b_n/q_{\lceil n/m\rceil}$
is comparable to $\log(en)$ and tends to infinity. Hence
$\I(Q)\subsetneq\I(A_Q)$.  Nevertheless
$Q$ and $A_Q$ belong to $\mathcal S_p$ for exactly the same
values of $p>0$, namely all $p\ge1$.
Thus this application detects algebraic ideal information invisible
to every Schatten-membership test.
\end{example}

\section{Why the affine coset is needed}\label{info:section}

\begin{proposition}[One unshifted ideal, no Borel complete coding]
\label{info:coset-essential}
Fix positive trace-class $Q$ of infinite rank, a rank-two skew kernel
$K=f\otimes g-g\otimes f$ with $f,g$ orthonormal, and
$0<\varepsilon<1$. There are centered Gaussian curvature arrays
$Z_\theta$, indexed by $\theta\in\{0,1\}^{\N}$, with standard
Gaussian row marginals and $\|\Gamma_\theta-I\|\le\varepsilon$,
whose raw and centered evaluations satisfy
\[
 \I(W_\theta)=\I(Q),\qquad
 \I(A_\theta+B)=\mathcal F+\I(B)
 \quad\text{a.s. for each fixed compact }B.
\]
Their complete centered responses coincide exactly when their
indices differ finitely. This equivalence admits no Borel complete
invariant in a standard Borel space.
\end{proposition}
\begin{proof}
Work on the support of $Q$, put $a_i=\sqrt{q_i}$ and
$J_0=\left(\begin{smallmatrix}0&1\\-1&0\end{smallmatrix}\right)$,
and set
\[
 c_i(\theta)=\tfrac\varepsilon2(1+\theta_i),\qquad
 S_\theta=\bigoplus_{i\ge1}c_i(\theta)J_0,
 \qquad\Gamma_\theta=I+S_\theta\otimes K.
\]
The perturbation is self-adjoint of norm at most $\varepsilon$
and has zero diagonal spatial blocks, proving the covariance and
row-marginal claims. Direct contraction gives
\[
 M_\theta=\E A_\theta
 =2D_aS_\theta D_a
 =2\bigoplus_{i\ge1}c_i(\theta)
       \sqrt{q_{2i-1}q_{2i}}\,J_0.
\]
Every $A_\theta$ has rank at most two and an atomless law:
its $(1,2)$ entry is a nonconstant quadratic polynomial of a
nondegenerate Gaussian vector. Finite-rank perturbation therefore gives
\[
 \I(W_\theta+B)=\mathcal F+\I(B-M_\theta),\qquad
 \I(A_\theta+B)=\mathcal F+\I(B)\quad\text{a.s.}
\]
The inequalities $q_{2i}\le\sqrt{q_{2i-1}q_{2i}}\le q_{2i-1}$
show by finite dilation that the paired weights generate $\I(q)$, and
$\varepsilon/2\le c_i(\theta)\le\varepsilon$, so
$\I(M_\theta)=\I(q)$. Two responses agree precisely when
$M_\theta-M_\eta$ has finite rank, equivalently when
$\theta,\eta$ differ at finitely many indices.

If a Borel complete invariant $h$ existed, each Borel inverse image
would be a tail event under fair Bernoulli measure, hence have probability
zero or one. A countable separating family in the target makes $h$
almost surely constant. Its full-measure fiber would then be one
countable finite-change class, which has measure zero: a contradiction.
\end{proof}
Zero Bochner means therefore do not force zero minimizing cosets.
The continuous covariance map still selects compact representatives;
equality of responses compares those representatives modulo the
fluctuation ideal.

\section*{Declaration of AI Use}
The author used GPT interactively in developing and
refining the mathematical arguments and the presentation of this paper.
The author takes full responsibility for the correctness and mathematical
content of the manuscript.

\end{document}